\documentclass[10pt]{amsart}
\usepackage{amsmath}
\usepackage{amssymb}
\usepackage{amsfonts}
\usepackage{amsthm}
\usepackage{enumerate}
\usepackage[all]{xy}
\usepackage[latin1]{inputenc}
\usepackage{graphicx}
\usepackage{latexsym}
\usepackage{accents}
\usepackage{mathrsfs}
\usepackage{indentfirst}
\usepackage{fancyhdr}
\input xy

\makeatletter
    \newtheoremstyle{definition}
        {5pt}
        {3pt}
        {}
        {0pt}
        {\scshape}
        {.}
        {5pt}
        {\thmname{#1} \thmnumber{#2} \thmnote{[#3]}} 

\newtheoremstyle{theorems}
        {5pt}
        {3pt}
        {\itshape}
        {0pt}
        {\scshape}
        {.}
        {5pt}
        {\thmname{#1} \thmnumber{#2}\thmnote{[#3]}} 

\swapnumbers

\theoremstyle{theorems}

\newtheorem{Theo}{Theorem}[section]
\newtheorem{Prop}[Theo]{Proposition}
\newtheorem{Cor}[Theo]{Corollary}
\newtheorem{Lemma}[Theo]{Lemma}
\newtheorem{Prop(BG)}[Theo]{Proposition (Bongartz-Gabriel)}
\newtheorem{Lemma(Asashiba)}[Theo]{Lemma(Asashiba)}
\newtheorem{Lemma(Gab)}[Theo]{Lemma(Gabriel)}
\newtheorem{Theo(Aus)}[Theo]{Theorem (Auslander)}
\newtheorem{Prop(RVDB)}[Theo]{Proposition (Reiten, van der Bergh)}

\theoremstyle{definition}
\newtheorem{Defn}[Theo]{Definition}

\newtheorem{Defn(Asashiba)}[Theo]{Definition (Asashiba)}

\newcommand{\Z}{\mathbb{Z}}

\newcommand{\La}{\mit\Lambda}
\newcommand{\cA}{\mathcal{A}}
\newcommand{\cB}{\mathcal{B}}
\newcommand{\cC}{\mathcal{C}}

\newcommand{\sC}{{\hspace{.5pt}\mathcal{C}}}
\newcommand{\oC}{\hspace{1.5pt}\overline{\hspace{-1pt}\mathcal{C}\hspace{-.2pt}}}
\newcommand{\uC}{\underline{\hspace{.6pt}\mathcal{C}\hspace{-1.5pt}}\hspace{1.5pt}}
\newcommand{\bC}{{\hspace{1pt}\mathcal{C}}}

\newcommand{\Hom}{{\rm Hom}}
\newcommand{\Ext}{{\rm Ext}}
\newcommand{\End}{{\rm End}}

\newcommand{\uHom}{\hspace{.8pt}\underline{\hspace{-.8pt}\rm Hom\hspace{-.6pt}}\hspace{.6pt}}
\newcommand{\oHom}{\overline{\rm Hom\hspace{-.8pt}}\hspace{.6pt}}

\def\Mod{\hbox{{\rm Mod}\hspace{0.5pt}}}

\def\ModLa{\hbox{{\rm Mod}{\hskip 0.5pt}\La}}

\def\mf{\mathfrak}

\newcommand{\id}{{\rm 1\hspace*{-0.56ex}\rule{0.1ex}{1.52ex}\hspace*{0.2ex}}}
\def\Ga{\hbox{$\it\Gamma$}}

\def\Sa{\hbox{${\it\Sigma}$}}
\def\Oa{\hbox{${\it\Omega}$}}
\def\La{\hbox{$\it\Lambda$}}
\def\Ta{\hbox{${\it\Theta}$}}

\newcommand{\dt}{{\accentset{\hspace{1pt}\mbox{\large\bfseries .}}{}}}
\newcommand{\wdt}{{\accentset{\hspace{4pt}\mbox{\small\bfseries .}}{}}}
\newcommand{\cdt}{\dt\hspace{2pt}}
\newcommand{\sdt}{\wdt\hspace{2pt}}

\newcommand{\pdt}{{\hspace{1.3pt}\dt\hspace{1pt}}}
\newcommand{\ydt}{{\hspace{0.5pt}\dt\hspace{1.5pt}}}

\begin{document}

\title[Almost split sequences]{\sc Almost split sequences in tri-exact categories}

\author[Shiping Liu]{Shiping Liu \vspace{-15pt}}

\address{Shiping Liu\\ D\'epartement de math\'ematiques, Universit\'e de Sherbrooke, Sherbrooke, Qu\'ebec, Canada}
\email{shiping.liu@usherbrooke.ca}

\author[Hongwei Niu]{\hspace{2pt} Hongwei Niu}

\address{Hongwei Niu\\D\'epartement de math\'ematiques, Universit\'e de Sherbrooke, Sherbrooke, Qu\'ebec, Canada}
\email{hongwei.niu@usherbrooke.ca}



\thanks{The first named author is supported in part by the Natu\-ral Science and Engineering Research Council of Canada.}

\maketitle

\begin{abstract}

We shall study the existence of almost split sequences in tri-exact categories, that is, extension-closed subcategories of triangulated cate\-gories. Our results unify and extend the existence theorems for almost split sequences in abelian categories and exact categories (that is, extension-closed subcategories of abelian categories), and those for almost split triangles in triangulated cate\-gories in \cite{AUS,Hap1,Kra2,LeZ,LNP,RvdB}. As applications, we shall obtain some new results on the existence of almost split sequences in the derived categories of all modules over an algebra with a unity or a locally finite dimensional algebra given by a quiver with relations. 

\vspace{5pt}

\end{abstract}

\section*{Introduction}

\medskip

Since its introduction in the last seventies; see \cite{AuR1, AuR2}, the Auslander-Reiten theory of almost split sequences has been playing a fundamental role in the modern representation theory of algebras; see, for example, \cite{ASS,ARS}. Later, Happel introduced the analogous theory of almost split triangles in triangulated categories; see \cite{Hap2, H}, making the Auslander-Reiten theory applicable in other areas of mathematics such as algebraic topology and algebraic geometry; see \cite{Rei,JOR,JOR2}. Since then, this theory has been further developed separately for exact categories 
and triangulated categories; 
see \cite{AUS,Hap1,Kra2, KLe, LeZ,LNP,RvdB}. Our purpose is to unify and extend these results by working with tri-exact categories. Observe that the existence of almost split sequences in a Krull-Schmidt category will help us to classify the indecomposable objects and describe certain morphisms in terms of the Auslander-Reiten quiver; see \cite{Liu1}. We shall outline the content of the paper section by section.

\medskip

In Section 1, in addition to laying down the foundation, we shall also study modules over an $R$-algebra, which are reflexive with respect to the minimal injective co-generator for ${\rm Mod}\hspace{.4pt}R$, where $R$ is a commutative ring. These modules will play the same role as those of finite length over an artin algebra. In case the algebra is reflexive and noetherian, we shall establish a duality between the $R$-noetherian modules and the $R$-artinian modules; see (\ref{Noe-Ref}), which generalizes the well-known Matlis duality; see \cite{AUS, BH}. 

\medskip

In Section 2, we shall study mainly the stable categories of a tri-exact category. The stable categories were first considered by Auslander and Reiten for modules over an artin algebra in order to establish the existence of almost split sequences; see \cite{AuR2}. Later, Lenzing and Zuazua defined the stable categories of an abelian category without projective or injection objects; see \cite{LeZ}, which carry over easily to an exact category; see \cite{LNP}. We shall extend them to tri-exact categories and show that every exact category is equivalent to a tri-exact category with equivalent stable categories. This ensures that the study of almost split sequences in exact categories and abelian categories is covered under our tri-exact setting. We should point out that a triangulated category coincides with its stable categories.

\medskip

In Section 3, we shall study the existence of an individual almost split sequence in a tri-exact category. Historically, one derives an almost split sequence from an Auslander-Reiten formula in an abelian category; see \cite{AUS} and \cite[(1.1)]{LeZ}, 
and an almost split triangle from a Serre formula in a triangulated category; see \cite[(2.2)]{Kra2} and \cite[(I.2.3)]{RvdB},
but the converses do not hold in general. These formulae involve taking the ``dual" of some stable Hom-spaces against injective modules over various rings; see \cite{AUS,ARS,Kra2,LeZ}. Recently, some necessary and sufficient conditions were found for the existence of an almost split sequence in an exact $R$-category, where the ``dual" is taken against an injective co-generator for ${\rm Mod}\hspace{.4pt}R$; see \cite[(2.2)]{LNP}. By taking the ``dual" against injective modules over rings mapping to the stable endomorphisms of two prescribed objects, we shall obtain some necessary and sufficient conditions for the existence of an almost split sequence in a tri-exact category; see (\ref{existence}), which essentially cover all the previously mentioned results.

\medskip

In Section 4, we shall be concerned with the global existence of almost split sequences in a tri-exact category. It is known that an Ext-finite abelian $R$-category with $R$ being artinian has almost split sequences if and only if it admits an Auslander-Reiten duality; see \cite{ARS,GR,LeZ} and a Hom-finite triangulated category over a field has almost split triangles on the right (or left) if and only if it admits a right (or left) Serre functor; see \cite[(I.2.3)]{RvdB}. We shall deal this problem for Hom-reflexive Krull-Schmidt tri-exact categories. This class of categories includes the category of noetherian modules and that of artinian modules over a noetheiran $R$-algebra with $R$ being notherian complete local, which are not Hom-finite if the algebra is not artinian; see \cite{AUS}. We shall show that such a tri-exact $R$-category has almost split sequences on the right (or left) if and only if it admits a full right (or left) Auslander-Reiten functor; see (\ref{ARS-ARF}). In the right (or left) triangulated case, the existence of almost split sequences on the right (or left) is equivalent to the existence of a right (or left) Auslander-Reiten functor, or equivalently, a right (or left) Serre functor with a proper image; see (\ref{AR-Serre}).

\medskip

In Section 5, we shall study the existence of almost split triangles in the derived categories of an abelian category with enough projective objects and enough injective objects. This has been done for the bounded derived category of finite dimensional modules over a finite dimensional algebra; see \cite{H, Hap2}. 
In the most general case, we shall show that an almost split triangle in the bounded derived cate\-gory starts with a bounded complex of injective objects and ends with a bounded complex of projective objects; see (\ref{ART-nec}) and (\ref{ART-bd}). In case the abelain category admits a Nakayama functor with respect to a subcategory of projective objects; see (\ref{Naka-Func}), we shall establish an existence theorem of an almost split triangle in the bounded derived category; 
see (\ref{ART-general}). In case the subcategory of projective objects is Hom-reflexive, we shall describe all possible almost split triangles in the bounded derived category; see (\ref{ART-refl}). These results will be applicable to the derived categories of module over an general algebra, a reflexive notherian algebra, or a locally finite dimensional algebra given by a quiver with relations.

\smallskip

\section{Preliminaries}

\medskip

\noindent The main objective of this section is to fix the notation and the terminology, which will be used throughout this paper, and collect some prelimi\-nary results. However, we shall also obtain some new results on modules over an algebra. Throughout this paper, morphisms in any category are composed from the right to the left.

\medskip

\noindent{\sc 1) Modules.} All rings and algebras except for those given by a quiver with relations have an identity. Let $\Sa$ be a ring or an algebra. We shall denote by $\Mod \Sa$ the category of all left $\Sa$-modules, and by ${\rm mod}\Sa$ the full subcategory of $\Mod \Sa$ of modules of finite length. For convenience, we shall identify the category of all right $\Sa$-modules with $\Mod \Sa^{\rm op}$, where $\Sa^{\rm op}$ is the opposite ring or the opposite algebra of $\Sa$.
A map $f: M\to N$ in $\Mod \Sa$ is called {\it socle essential} provided that ${\rm Im}(f) \cap {\rm Soc}(N)$ is non-zero whenever ${\rm Soc}(N)$ is non-zero.

\medskip

Let $M$ be a left or right $\Sa$-module. Then $M^*=\Hom_{\it\Sigma}(M, \Sa)$ is a right or left $\Sa$-module, respectively. Given $u\in M$, we have
$\hat{u}\in M^{**}=\Hom_{\it\Sigma}(M^*, \Sa)$, sending $f\in M^*$ to $f(u)$. The map $\rho_{_M}: M\to M^{**}$, sending $u$ to $\hat{u},$ is clearly a natural $\Sa$-linear map. It is well known; see, for example, \cite[(3.15)]{ROT} that $M$ is finitely generated projective if and only if it has a finite  projective basis $\{u_i; f_i\}_{1\le i\le n}$, where $u_i\in M$ and $f_i\in M^*$, such that $u=\sum_{i=1}^n\, f_i(u) u_i$ (or $u=\sum_{i=1}^n\, u_i f_i(u))$ for all $u\in M$. We shall denote by ${\rm proj}\hspace{.4pt}\Sa$ the full subcategory of $\Mod \Sa$ of finitely generated projective modules. The following statement is probably well-known.

\begin{Lemma}\label{fgp-equiv}

Let $\Sa$ be a ring or an algebra. Then $\Hom_{\it\Sigma}(-, \Sa): {\rm proj}\hspace{.5pt}\Sa \to {\rm proj} \hspace{.5pt}\Sa^{\rm op}$ is a duality.

\end{Lemma}

\noindent{\it Proof.} Let $P\in {\rm proj}\hspace{.5pt}\Sa$ with a finite projective basis $\{u_i; f_i\}_{1\le i\le n}$. Given $f \in P^*$ and $u\in P,$ we obtain
$$
\textstyle(\sum_{i=1}^n f_i \, \hat{u}_i(f))(u)
= \sum_{i=1}^n \left(f_i f(u_i)\right)(u) = \sum_{i=1}^n f_i(u) f (u_i) = f (\sum_{i=1}^n f _i(u)u_i).$$
Since $u=\sum_{i=1}^n f _i(u)u_i$, we conclude that $f=\sum_{i=1}^n f_i \, \hat{u}_i(f)$. That is, $\{f_i; \hat{u}_i\}_{1\le i\le n}$ is a projective basis of $P^*$. In particular, $P^*\in {\rm proj}\hspace{.5pt}\Sa^{\rm op}$. If $u\in P$ is non-zero, since $u=\sum_{i=1}^n f_i(u) u_i,$ we see that $\hat{u}$ is non-zero. That is, $\rho_{_P}$ is a monomorphism. As shown above, $P^{* *}$ has a projective basis $\{\hat{u}_i; \hat{f}_i\}_{1\le i\le n}$. Given $\varphi\in P^{**}$, we obtain
$$\textstyle \varphi=\sum_{i=1}^n \hat{f}_i(\varphi) \hat{u}_i=\sum_{i=1}^n \varphi(f_i) \hat{u}_i=\rho_{_P}(\sum_{i=1}^n \varphi(f_i) u_i).$$ Thus, $\rho_{_P}$ is an isomorphism.
The proof of the lemma is completed.

\medskip

Throughout this paper, $R$ will stand for a commutative ring and $I_{\hspace{-1pt}R}$ for a minimal injective co-generator for $\Mod R$; see \cite[(18.19)]{AnF}. We shall use frequently the functor $D=\Hom_R(-, I_{\hspace{-1pt}R}): \Mod R \to \Mod R$. The following statement is probably known.


\begin{Prop}\label{Mod-Commu}

Let $U$ be a module over a commutative ring $R$.

\vspace{-1pt}

\begin{enumerate}[$(1)$]

\item If $U$ is of finite length $n$, then $DU$ is also of length $n$.

\item If $U$ is finitely co-generated, then $DU$ is finitely generated.

\item If $DU$ is artinian or noetherian, then $U$ is noetherian or artinian respectively.


\end{enumerate}

\end{Prop}

\noindent{\it Proof.} Assume first that $U$ is simple. In particular, $U=Ru$ for some $u\in U$. Consider some non-zero linear functions $f, g\in DU$. Since $I_{\hspace{-1pt}R}$ is a minimal injective co-generator, ${\rm soc}(I_{\hspace{-1pt}R})$ contains exactly one copy of $U$; see \cite[(18.19)]{AnF}. Therefore, $g(U)=f(U)$, and hence, $g(u)= r f(u)$ for some $r\in R$. This yields $g=rf$. Thus, $DU$ is also simple. By induction, we can establish Statement (1).

Assume next that $U$ is finitely co-generated, that is, $U$ has an essential socle $S=S_1\oplus \cdots\oplus S_t$, where the $S_i$ are simple; see \cite[(10.4)]{AnF}. Consider the canonical projections $p_i: S\to S_i$ and the canonical injections $q_i: S_i\to S$, and fix some monomorphisms $f_i: S_i \to I_{\hspace{-1pt}R}$, for $i=1, \ldots, t$. Letting $q: S\to U$ be the inclusion, we obtain $R$-linear maps $g_i: U\to I_{\hspace{-1pt}R}$ such that $g_i q =f_i p_i$, for $i=1, \ldots, t$. Given any $R$-linear map $g: U\to I_{\hspace{-1pt}R}$, as seen above, $gqq_i=r_i f_i$ for some $r_i\in R$. This yields $gq=\textstyle\sum_{i=1}^t gqq_i p_i=\sum_{i=1}^t r_if_ip_i=( \sum_{i=1}^t r_i g_i)q.$ Since $q$ is an essential monomorphism, $g=\sum_{i=1}^t r_i g_i.$  Statement (2) is established.

Finally, given a submodule $V$ of $U$, we denote by $V^\perp$ the submodule of $DU$ of $R$-linear maps vanishing on $V$ and by $^\perp(V^\perp)$ the submodule of $U$ of elements annihilated by the $R$-linear maps in $V^\perp$. Then, $V\subseteq {}^\perp(V^\perp)$. We claim that $V= {}^\perp(V^\perp)$. Otherwise, we can find an $R$-linear map $h: U\to I_{\hspace{-1pt}R}$ such that $h({}^\perp(V^\perp))\ne 0$ but $h(V)=0$, contrary to the definition. Using this claim, we may easily establish Statement (3). The proof of the proposition is completed.

\medskip

Let $A$ be an $R$-algebra. A left or right $A$-module $M$ is called {\it $R$-noetherian} or {\it $R$-artinian} if $_RM$ is noetherian or artinian; and $A$ is called a {\it noetherian} or {\it reflexive $R$-algebra} if $_AA$ is $R$-noetherian or $R$-reflexive, respectively. Note that our definition of a noetherian $R$-algebra is different from the classical one, where $R$ is assumed to be noetherian. Consider the exact functors $D=\Hom_R(-, I\hspace{-1.5pt}_R)\hspace{-1pt}:\hspace{-1pt} \Mod A \hspace{-1pt}\to\hspace{-1pt} \Mod A^{\rm op}$ and $D=\Hom_R(-, I\hspace{-1.5pt}_R)\hspace{-1pt}:\hspace{-1pt} \Mod A^{\rm op} \hspace{-1pt}\to\hspace{-1pt} \Mod A$. Given a left or right $A$-module $M$, we obtain a canonical $A$-linear monomorphism $\sigma_{\hspace{-1.5pt}_M}: M\to D^2M$ so that $\sigma_{\hspace{-1pt}_M}(x)(f)=f(x)$, for $x\in M$ and $f\in DM$. We shall say that $M$ is {\it $R$-reflexive} if $\sigma_{\hspace{-1pt}_M}$ is bijective.
%
%

\begin{Lemma}\label{alg-ref-mod}

Let $A$ be an $R$-algebra. The full subcategory ${\rm RMod}\hspace{.4pt}A$ of $\Mod A$ of $R$-reflexive modules is abelian, contains all modules of finite $R$-length, and admits a duality $D: {\rm RMod}\hspace{.5pt}A \to {\rm RMod}\hspace{.5pt}A^{\rm op}$.

\end{Lemma}

\noindent{\it Proof.} Considering the canonical monomorphisms and applying the Snake Lemma, we see that  ${\rm RMod}\hspace{.5pt}A$
is closed under taking submodules and quotient modules, that is, it is abelian. If $M\in \Mod A$ is of $R$-length $n$, by Proposition \ref{Mod-Commu}(1), so is $D^2M$, and consequently, $\sigma_{\hspace{-1.5pt}_M}: M\to D^2M$ is an isomorphism. 
Finally, by the definition of reflexive modules, $D: {\rm RMod}\hspace{.5pt}A \to {\rm RMod}\hspace{.5pt}A^{\rm op}$ and $D: {\rm RMod}\hspace{.5pt}A^{\rm op} \to {\rm RMod}\hspace{.5pt}A$ are mutual quasi-inverse. The proof of the lemma is completed.

\medskip

Consider now the endofunctors $\nu_{\hspace{-1.5pt}_A}=D\Hom_A(-, A)$ and $\nu_{\hspace{-1.5pt}_A}^{\mbox{\hspace{1pt}-\hspace{-3pt}}}=\Hom_A(D(-), A)$ of $\Mod A$. Put ${\rm inj}\hspace{.4pt} A=\nu_{\hspace{-1.5pt}_A}({\rm proj}\hspace{.5pt}A)$ which, by Lemma \ref{fgp-equiv}, contains only injective mo\-dules. Let ${\rm mod}^{+\hspace{-3pt}}A$ stand for the full subcategory of $\Mod A$ of finitely generated modules, and ${\rm mod}^{-\hspace{-3.5pt}}A$ for that of modules finitely co-generated by ${\rm inj}\hspace{.4pt} A$.


\begin{Theo}\label{Noe-Ref} Let $A$ be a reflexive noetherian $R$-algebra.

\vspace{-3pt}

\begin{enumerate}[$(1)$]

\item The functors $\nu_{\hspace{-1.5pt}_A}: {\rm proj}\hspace{.5pt}A \to {\rm inj}\hspace{.5pt}A$ and $\nu_{\hspace{-1.5pt}_A}^{\mbox{\hspace{.5pt}-\hspace{-3pt}}}: {\rm inj}\hspace{.5pt}A\to {\rm proj}\hspace{.5pt}A$ are mutual quasi-inverses, where ${\rm proj}\hspace{.5pt}A$ and ${\rm inj}\hspace{.5pt}A$ have as objects all $R$-noetherian projective modules and all $R$-artinian injective modules, respectively.

\vspace{.5pt}

\item There exists a duality $D=\Hom_R(-, I_{\hspace{-1pt}R}): {\rm mod}^{+\hspace{-3pt}}A^{\rm op}\to {\rm mod}^{-\hspace{-3pt}}A$, where
${\rm mod}^{+\hspace{-3pt}}A$ and ${\rm mod}^{-\hspace{-3.5pt}}A$ are abelian subcategories of ${\rm RMod}\hspace{.4pt}A$, whose objects are all $R$-noetherian modules and all $R$-artinian modules, respectively.


%
%
%

\end{enumerate}

\end{Theo}

\noindent{\it Proof.} Since $_AA$ is $R$-reflexive and $R$-noetherian, we deduce from Lemma \ref{alg-ref-mod} that ${\rm mod}^{+\hspace{-3pt}}A$ is an abelian subcategory of ${\rm RMod}\hspace{.5pt}A$, whose objects are clearly the $R$-noetherian $A$-modules. In particular, the objects of ${\rm proj}\hspace{.5pt}A$ are the $R$-noetherian projective $A$-modules.
On the other hand, since $A^{\rm op}$ is also a reflexive noetherian $R$-algebra, ${\rm mod}^{+\hspace{-3pt}}A^{\rm op}$ is an abelian subcate\-gory of ${\rm RMod}\hspace{.5pt}A^{\rm op}$. Consider the equiva\-lence $\Hom_A(-, A): {\rm proj}\hspace{.4pt}A \to {\rm proj}\hspace{.4pt}A^{\rm op}$ and the duality $D: {\rm RMod}\hspace{.5pt}A^{\rm op} \to {\rm RMod}\hspace{.5pt}A$ in Lemmas \ref{fgp-equiv} and \ref{alg-ref-mod}, we see that ${\rm inj}\hspace{.5pt}A$ is a subcategory of ${\rm RMod}\hspace{.5pt}A$, whereas the functors $\nu_{\hspace{-2pt}_A}: {\rm proj}\hspace{.5pt}A\to {\rm inj}\hspace{.5pt}A$ and $\nu_{\hspace{-2pt}_A}^{\hspace{.5pt}\mbox{-\hspace{-3.5pt}}}: {\rm inj}\hspace{.5pt}A\to {\rm proj}\hspace{.5pt}A$ are mutual quasi-inverses.

Being finitely co-generated by ${\rm inj}\hspace{.5pt}A$, by Lemma \ref{alg-ref-mod}, ${\rm mod}^{-\hspace{-3pt}}A$ is a subcategory of ${\rm RMod}\hspace{.5pt}A$. Given $M\in {\rm RMod}\hspace{.4pt}A^{\rm op}$, in view of the duality $D: {\rm RMod}\hspace{.5pt}A^{\rm op}\to {\rm RMod}\hspace{.5pt}A$, we see that $M\in {\rm mod}^{-\hspace{-3pt}}A$ if and only if $DM\in {\rm mod}^{+\hspace{-3pt}}A^{\rm op}$. Thus, we obtain a duality $D: {\rm mod}^{+\hspace{-3pt}}A^{\rm op}\to {\rm mod}^{-\hspace{-3pt}}A$. In particular, ${\rm mod}^{-\hspace{-3pt}}A$ is abelian. Since ${\rm mod}^{+\hspace{-3pt}}A^{\rm op}$ contains only $R$-noetherian modules, by Proposition \ref{Mod-Commu}(3), ${\rm mod}^{-\hspace{-3pt}}A$ contains only $R$-artinian modules. On the other hand, if $M\in \Mod A$ is $R$-artinian, then $_RM$ is finitely co-generated, and by Proposition \ref{Mod-Commu}(2), $DM$ is finitely generated over $R$. In particular, $DM\in {\rm mod}^{+\hspace{-3pt}}A^{\rm op}$, and hence, $D^2M\in {\rm mod}^{-\hspace{-3pt}}A$. Since  ${\rm mod}^{-\hspace{-3pt}}A$ is abelian and $\sigma_{\hspace{-1pt}_M}: M\to D^2M$ is a monomorphism, $M\in {\rm mod}^{-\hspace{-3pt}}A$. Finally, let $I$ be an $R$-artinian injective $A$-module. Then, $I$ is an injective object in ${\rm mod}^{-\hspace{-3pt}}A$. Hence, $I\cong DP$, where $P$ is a projective object in ${\rm mod}^{+\hspace{-3pt}}A^{\rm op}$. It is then easy to see that $P\in {\rm proj}\hspace{.4pt}A^{\rm op},$ that is, $I\in {\rm inj}\hspace{.4pt}A.$ The proof of the theorem is completed.

\medskip

\noindent{\sc Remark.} A noetherian algebra over a commutative noetherian complete local ring is reflexive; see \cite[Section 5]{AUS}. Thus, Theorem \ref{Noe-Ref} generalizes the well-known Matlis duality; see \cite[Section 5]{AUS} and \cite[(3.2.13)]{BH}.

\medskip

We conclude this subsection with algebras given by a quiver with relations. Let $Q$ be a locally finite quiver with vertex set $Q_0$. An infinite path in $Q$ is called {\it left infinite} if it has no starting point and {\it right infinite} if it has no ending point. Given a field $k$, an ideal $J$ in the path algebra $kQ$ is called {\it weakly admissible} if it lies in the ideal generated by the paths of length two; and {\it locally admissible} if, for any $x\in Q_0$, there exists $n_x\in \Z$ for which $R$ contains all paths of length $\ge n_x$, starting or ending with $x$. Consider $\La=kQ/J$, where $J$ is weakly admissible, with a complete set of pairwise orthogonal primitive idempotents $\{e_x \mid x\in Q_0\}$. One calls $\La$ {\it locally finite dimensional} if $e_x\La e_y$ is finite dimensional for all $x, y\in Q_0$; and  {\it strongly locally finite dimensional} if $J$ is locally admissible; see \cite[Section 1(4)]{BHL}.
Assume that $\La$ is locally finite dimensional. We shall denote by $\Mod\La$ the category of all left $\La$-modules $M$ such that $M=\oplus_{x\in Q_0}\, e_xM$, and by ${\rm mod}^{\hspace{.5pt}b\hspace{-2.5pt}}\La$ the full subcategory of $\Mod\La$ of finite dimensional modules. Given $x\in Q_0$, one obtains a projective module $P_x=\La e_x$ and an injective module $I_x=D(\La^{\rm op}e_x)$; see \cite[Section 3]{BHL}. Let ${\rm proj}\,\La$ and ${\rm inj}\hspace{.4pt}\La$ be
the strictly additive subcategories of $\Mod\La$ generated by the $P_x$ with $x\in Q_0$ and by the $I_x$ with $x\in Q_0$, respectively. Observe that $\La$ is strongly locally finite dimensional if and only if all $P_x$ and $I_x$ with $x\in Q_0$ are finite dimensional.

\medskip

\noindent{\sc 2) Additive categories.} Throughout this paper, all functors between additive cate\-gories are additive. Let $\mathcal{A}$ be an additive category. A {\it strictly additive} subcategory of $\cA$ is a full subcategory which is closed under finite direct sums, direct summands and isomorphisms. An object in $\mathcal{A}$ is called {\it strongly indecomposable} if it has a local endomorphism ring. One says that $\cA$ is {\it Krull-Schmidt} if every non-zero object is a finite direct sum of strongly indecomposable objects; and in this case, $\cA/\mathcal{I}$ is Krull-Schmidt, for every ideal $\mathcal{I}$ in $\cA$.

\medskip

A morphism $f: X\to Y$ in $\cA$ is called {\it left minimal} if every morphism $h: Y\to Y$ such that $f=hf$ is an automorphism;
%
%
%
{\it left almost split} if $f$ is not a section such that every non-section morphism $g: X \to M$ factors through it; and {\it minimal left almost split} if it is left minimal and left almost split. In the dual situations, one says that $f$ is {\it right minimal}, {\it right almost split}, and {\it minimal right almost split}, respectively.

\medskip

\noindent{\sc 3) Stable categories of exact categories.} Let $\mathscr{C}$ be an exact category, that is an extension-closed subcategory of an abelian category $\mf{A}$; see \cite[Section 2]{LNP}.
Given $X, Y\in \mathscr{C}$, one writes $\Ext_\mathscr{\hspace{.5pt}C}^1(X, Y)=\Ext_\mathfrak{\hspace{.5pt}A}^1(X, Y)$. A morphism $f: X\to Y$ in $\mathscr{C}$ is called {\it injectively trivial} in $\mathscr{C}$ if the push-out map
$$\Ext_\mathscr{\hspace{.5pt}C}^1(Z, f): \, \Ext_\mathscr{\hspace{.5pt}C}^1(Z, X)\to \Ext_\mathscr{\hspace{.5pt}C}^1(Z, Y): \delta\mapsto f \cdot \delta$$ vanishes for all $Z\in \mathscr{C}$; and {\it projectively trivial} in $\mathscr{C}$ if the pull-up map
$$\Ext_\mathscr{\hspace{.5pt}C}^1(f, Z): \, \Ext_\mathscr{\hspace{.5pt}C}^1(Y, Z)\to \Ext_\mathscr{\hspace{.5pt}C}^1(X, Z): \zeta\mapsto \zeta \cdot f$$ vanishes for all $Z\in \mathscr{C}$. Denoting by $\mathcal{I}_{\hspace{.5pt}\mathscr{C}}$ the ideal of injectively trivial morphisms and by $\mathcal{P}_{\mathscr{\hspace{.5pt}C}}$ that of projectively trivial morphisms, one obtains the {\it injectively stable category} $\hspace{2pt}\overline{\hspace{-2pt}\mathscr{C}} =\mathscr{C}/\mathcal{I}_{\hspace{.5pt}\mathscr{C}}$ and the {\it projectively stable category} $\underline{\mathscr{C}\hspace{-2pt}}\hspace{1pt}=\mathscr{C}/\mathcal{P}_{\mathscr{\hspace{.5pt}C}}$ of $\mathscr{C}$; see \cite{LNP, LeZ}. In the sequel, we shall write as usual
$\Hom_{\hspace{2.5pt}\overline{\hspace{-2pt}\mathscr{C}}\hspace{.5pt}}(X, Y)=\overline{\rm Hom}\hspace{.5pt}_\mathscr{C}(X, Y)$ and $\Hom_{\hspace{2.5pt}\underline{\hspace{-1pt}\mathscr{C}\hspace{-1pt}}\hspace{1pt}}(X, Y)=\underline{\hspace{-.5pt}\rm Hom\hspace{-.5pt}}_{\hspace{1pt}\mathscr{C}}(X, Y)$, for all $X, Y\in \mathscr{C}$.

\medskip

\noindent{\sc 4) Triangulated categories.} Let $\cA$ be a triangulated category, whose translation functor will always be written as $[1].$ An exact triangle \vspace{-3pt}
$$\xymatrix{X \ar[r]^f &Y \ar[r]^g & Z \ar[r]^-\delta &X[1]} \vspace{-2pt}$$ in $\cA$ is called {\it almost split} if $f$ is minimal left almost split and $g$ is minimal right almost split; see \cite[(4.1)]{H}. In this case, $X$ and $Z$ will be called the {\it starting term} and the {\it ending term}, respectively.

\begin{Defn}

A full subcategory $\cC$ of a triangulated category $\cA$ is called {\it exten\-sion-closed} provided, for any
exact triangle $\xymatrixcolsep{18pt}\xymatrix{X \ar[r] &Y \ar[r] &Z \ar[r] &X[1]}$ in $\mathcal{A}$, that $Y \in \mathcal{C}$ whenever $X, Z\in \mathcal{C}.$ In this case, we shall simply call $\cC$ a {\it tri-exact category.}

\end{Defn}

\medskip

Observe that an extension-closed subcategory of a triangulated category is strictly additive, but it is not necessarily closed under the translation functor $[1]$.


\begin{Defn}\label{tri-sub}

An extension-closed subcate\-gory $\mathcal{C}$ of a triangulated category $\cA$ will be called

\begin{enumerate}[$(1)$]

\item {\it left triangulated} provided, for any exact triangle
$\xymatrixcolsep{18pt}\xymatrix{\hspace{-3pt} X \ar[r] &Y \ar[r] &Z \ar[r] &X[1]\hspace{-3pt}}$ in $\mathcal{A}$, that $X\in \cC$ whenever $Y, Z\in \cC$.

\item {\it right triangulated} provided, for any exact triangle
$\xymatrixcolsep{18pt}\xymatrix{\hspace{-3pt}X \ar[r] &Y \ar[r] &Z \ar[r] &X[1]\hspace{-3pt}}$ in $\mathcal{A}$, that $Z\in \cC$ whenever $X, Y\in \cC$.

%

\end{enumerate}

\end{Defn}

\medskip

Let $\cC$ be a full subcategory of $\cA$. For $n\in \Z$, we denote by $\cC[n]$ the full subcategory of $\cA$ generated by the objects $X[n]$ with $X\in \cC$. Clearly, $X\in \cC[n]$ if and only if $X[-n]\in \cC$. Turning exact triangles in $\cA$, we obtain the following observation.

\begin{Lemma}\label{left-tri}

An extension-closed subcategory $\cC$ of a triangulated category is left triangulated if and only if $\cC[-1]\subseteq \cC;$ and right triangulated if and only if $\cC[1]\subseteq \cC.$

\end{Lemma}

%

\medskip

Observe that a left or right triangulated subcategory of a triangulated category is a left or right triangulated category as defined in \cite[(1.1)]{ABM} and will be simply called a {\it left triangulated category} with a left translation $[-1]$ or a {\it right triangulated category} with a  right translation $[1]$, respectively.

%

\medskip

\noindent{\sc 5) Derived categories.} Let $\mathcal{A}$ be an additive category. We shall denote by $C(\mathcal{A})$ the {\it complex category} of $\cA$. The full subcate\-gories of $C(\mathcal{A})$ of bounded-above complexes, of bounded-below complexes, and of bounded complexes will be written as $C^-(\mathcal{A})$, $C^-(\mathcal{A})$ and $C^b(\mathcal{A})$, respectively. Given $*\in \{+, -, b, \emptyset\}$, let $K^*(\mathcal{A})$ stand for the {\it homotopy category}, that is the quotient of $C^*(\cA)$ modulo the null-homotopic morphisms, which is triangulated with a canonical projection functor $\mathbb{P}^*: C^*(\cA)\to K^*(\cA)$; and  $D^*(\mathcal{A})$ for the {\it derived category}, that is the localization of $K^*(\mathcal{A})$ with respect to the quasi-isomorphisms, which is also triangulated with a canonical localization functor $\mathbb{L}^*: K^*(\mathcal{A}) \to D^*(\mathcal{A})$. Given a morphism $f^\ydt$ in $C^*(\cA),$ we shall write $\bar{f}^\ydt= \mathbb{P}^*(f^\ydt)\in K^*(\cA)$ and $\tilde{f}^\ydt=\mathbb{L}^*(\bar{f}^\ydt)\in D^*(\cA).$ Given a complex $M^\ydt\in C^b(\mathcal{A}),$ its {\it width} $w(M^\ydt)$ is an integer defined by $w(M^\ydt)=0$ if $M^i=0$ for all $i\in \Z$; and otherwise, $w(M^\ydt)=t-s+1$, where $s\le t$ such that $M^s$ and $M^t$ are non-zero, but $M^i=0$ for all $i\notin [s, t]$.
Furthermore, given a complex $(X^\ydt, d^\sdt)$ over $\cA$ and some integer $n$, one defines two {\it brutal truncations} \vspace{-5pt}
$$\kappa_{\ge n}(X^\cdt): \; \textstyle \xymatrix{\cdots \ar[r] &0 \ar[r] &X^{n} \ar[r]^{d^n} &X^{n+1} \ar[r]^{d^{n+1}} &X^{n+2} \ar[r] &\cdots}\vspace{-3pt}$$
and \vspace{-10pt}
$$\kappa_{\le n}(X^\cdt): \; \xymatrix{\cdots \ar[r] &X^{n-2}\ar[r]^{d^{n-2}} &X^{n-1} \ar[r]^{d^{n-1}} &X^{n} \ar[r] &0 \ar[r] &\cdots,}$$ where $X^n$ is the component of degree $n$ in both complexes, with two canonical morphisms $\mu_n^\ydt: \kappa_{\ge n}(X^\ydt) \to X^\ydt \vspace{1pt}$ such that $\mu_n^p=1_{_{X^p}}$ for $p\ge n$ and $\mu_n^p=0$ for $p<n;$ and $\pi_n^\ydt: X^\ydt \to \kappa_{\le n}(X^\ydt)$ such that $\pi_n^p=1_{_{X^p}}$ for $p\le n$ and $\pi_n^p=0$ for $p > n.$ The following statement is well-known; see \cite[(III.4.4.2)]{Mil} and \cite[(1.3)]{Hap1}.


\begin{Lemma}\label{truncation-exact}

Let $\mathcal{A}$ be an additive category. \vspace{-2pt} If $X^\ydt\in C(\cA)$ and $n\in \Z,$ then $K(\mathcal{A})$ has an exact triangle
$\xymatrixcolsep{18pt}\xymatrix{\hspace{-3pt}\kappa_{\ge n}(X^\cdt) \ar[r] 
& X^\cdt \ar[r] 
& \kappa_{\le n-1}(X^\cdt) \ar[r] 
&\kappa_{\ge n}(X^\cdt)[1].}$

\end{Lemma}

\smallskip

Consider now the derived category of an abelian category $\mf{A}$. Fix an integer $n$. We shall denote by $D^{\le n}(\mf A)$ and $D^{\ge n}(\mf A)$ the full subcategories of $D(\mf A)$ of complexes $M^\ydt$ with ${\rm H}^i(M^\ydt)=0$ for all $i> n$ and of complexes $M^\ydt$ with ${\rm H}^i(M^\ydt)=0$ for all $i< n$ respectively, where ${\rm H}^i(M^\ydt)$ is the $i$-th cohomology of $M^\ydt$. By Lemma \ref{tri-sub}, $D^{\le n}(\mf A)$ is right triangulated and $D^{\ge n}(\mf A)$ is left triangulated in $D(\mf A)$. Let $(X^\pdt, d^\pdt)$ be a complex over $\mf{A}$. Writing $d^n=q^n p^n$, where $p^{n}: X^{n} \to C^{n}$ is the cokernel of $d^{n-1}$, and $d^{n-1}=i^nj^{n-1}$, where $i^n: K^n\to X^n$ is the kernel of $d_n$, we obtain two {\it smart truncations} \vspace{-3pt}
$$\tau_{\ge n}(X^\cdt): \quad \xymatrix{\cdots \ar[r] &0\ar[r]& C^n \ar[r]^{q^n} &X^{n+1} \ar[r]^{d^{n+1}} &X^{n+2} \ar[r] &\cdots,}
$$
where $C^n$ is of degree $n$, with a canonical projection $p_n^\sdt: X^\ydt\to \tau_{\ge n}(X^\ydt)$ so that $p_n^n= p^n$ and $p^s_n=\id_{X^s}$ for all $s>n$; and \vspace{-3pt}
$$\tau_{\le n}(X^\cdt): \quad \xymatrix{\cdots \ar[r] & X^{n-2} \ar[r]^{d^{n-2}} & X^{n-1} \ar[r]^-{j^{n-1}} & K^n \ar[r] & 0\ar[r] & \cdots }$$
where $K^n$ is of degree $n$, with a canonical injection $i_n^\sdt: \tau_{\le n}(X^\ydt)\to X^\ydt$ so that $i_n^n=i^n$ and $i_n^t=\id_{X^t}$ for all $t<n.$



%
%
%

%


\begin{Lemma}\label{bhomology}

Let $\mf B$ be an abelian subcategory of an abelian category $\mf{A}$. Consider $X^\ydt \in C^*(\mf A)$ with $*\in \{\emptyset, -, +, b\}$ and $Y^\ydt\in C(\mf B)$. If $X^\ydt\cong Y^\ydt$ in $D(\mf A)$, then there exists $Z^\ydt\in C^*(\mf B)$ such that $X^\ydt\cong Z^\ydt$ in $D^*(\mf A)$.

\end{Lemma}

\noindent{\it Proof.} We shall only consider the case where $X^\ydt \in C^b(\mf A)$ and $Y^\ydt\in C(\mf B)$ such that $X^\ydt\cong Y^\ydt$ in $D(\mf A)$. Let $s, t$ with $s\le t$ be such that $X^p=0$ for $p\not\in [s, t]$. Then, ${\rm H}^p(Y^\ydt)=0$ for all $p\not\in [s, t]$, and hence, the canonical injection
$i_t^\sdt: \tau_{\le t}(Y^\ydt)\to Y^\ydt$ and the canonical projection $p_s^\sdt: \tau_{\le t}(Y^\ydt) \to \tau_{\ge s}(\tau_{\le t}(Y^\ydt))$ are quasi-isomorphisms; see \cite[(III.3.4.1), (III.3.4.2)]{Mil}. As a consequence, $\tau_{\ge s}(\tau_{\le t}(Y^\ydt))\in C^b(\mf B)$ such that $X^\cdt\cong \tau_{\ge s}(\tau_{\le t}(Y^\ydt))$ in $D^b(\mf A)$. The proof of the lemma is completed.

\medskip

In the same fashion, one can show that $D^*(\mf A)$ with $*\in \{-, +, b\}$ fully embeds in $D(\mf A)$; see \cite[(III.3.4.3), (III.3.4.4), (III.3.4.5)]{Mil}. In the sequel, we shall always regard $D^*(\mf A)$ as a full triangulated subcategory of $D(\mf A)$.


\section{Tri-exact structure and stable categories}

\medskip

\noindent The objective of this section is to study the tri-exact structure in a non-axiomatic fashion and introduce the stable categories of a tri-exact category. They are analo\-gous to the exact structure and the stable categories of an exact category as described in \cite{LeZ, LNP}. More importantly, we shall show that every exact category is equivalent to a tri-exact category with equivalent stable categories.

\medskip

Throughout this section $\mathcal{C}$ will denote a tri-exact category, say an extension-closed subcategory of a triangulated category $\cA$. Given $X, Y\in \cC$, we shall write $\Ext^1_\mathcal{C}(X, Y)=\Hom_\mathcal{A}(X, Y[1])$, whose elements will be called {\it extensions} of $Y$ by $X$. Given an extension $\delta\in{\rm Ext}^1_\mathcal{C}(X, Y)$ and two morphisms $f\in \Hom_{\cC}(M, X)$ and $g\in \Hom_{\cC}(Y, N)$, we shall define $$\delta \cdot f=\delta \circ f \in \Hom_\cA(M, Y[1])=\Ext_\cC^1(M, Y)\vspace{-6pt}$$ and \vspace{-2pt}
$$g\cdot \delta=g[1]\circ \delta\in \Hom_\cA(X, N[1])=\Ext_\cC^1(X, N). \vspace{4pt}$$

This yields the following trivial observation.

\begin{Lemma}\label{Ext-bimod}

Let $\cC$ be a tri-exact category. Given an extension $\delta$ and two morphisms $f, g$ in $\cC$, the following equations hold whenever the composites make sense$\,:$ \vspace{-8pt}
$$(g\cdot \delta )\cdot f =g\cdot (\delta \cdot f); \quad (\delta\cdot f) \cdot g  =\delta \cdot (fg); \quad g \cdot (f \cdot \delta)= (gf ) \cdot \delta.$$

\end{Lemma}

\medskip

The tri-exact structure of $\cC$ consists of the tri-exact sequences as defined below.


\begin{Defn}

Let $\cC$ \vspace{-.5pt} be an extension-closed subcategory of a triangulated cate\-gory $\cA$.  A sequence of morphisms  $\hspace{-3pt}\xymatrixcolsep{18pt}\xymatrix{X \ar[r] & Y \ar[r] & Z}\hspace{-2pt}$ in $\cC$ is called a {\it tri-exact sequence} if it embeds in an exact triangle \vspace{-2pt}
$$\xymatrix{X \ar[r] &Y \ar[r] &Z \ar[r]^-\delta & X[1]}\vspace{-5pt}$$ in $\cA$. In this case, we say that the tri-exact sequence is {\it defined} by $\delta\in \Ext_\cC^1(Z, X)$, and call $X$ the {\it starting term} and $Z$ the {\it ending term}.

\end{Defn}

\medskip

\noindent{\sc Remark.} It is well-known; see \cite[(1.2)]{H} that a tri-exact sequence is a pseudo-exact sequence as defined in \cite[Section 1]{Liu1}. The converse, however, is not true.

\medskip

The following statement describes some basic properties of the tri-exact structure of a tri-exact category.

\begin{Lemma}\label{section-retraction}

Let $\cC$ be a tri-exact category, and let $\hspace{-2pt}\xymatrixcolsep{18pt}\xymatrix{X\ar[r]^f &Y\ar[r]^g &Z}\hspace{-3pt}$ be a tri-exact sequence defined by an extension $\delta\in \Ext_\mathcal{C}^1(Z, X)$.

\vspace{-1pt}

\begin{enumerate}[$(1)$]

\item A morphism $u: X\to M$ factors through $f$ if and only if $u\cdot \delta =0$.

\item A morphism $v: N\to Z$ factors through $g$ if and only if $\delta \cdot v=0.$

\item The morphism $f$ is a section if and only if $g$ is a retraction if and only if $\delta=0$.

\item If $X$ or $Z$ is strongly indecomposable, then $g$ is right minimal or $f$ is left minimal, respectively.

\end{enumerate}

\end{Lemma}

\noindent{\it Proof.} Assume that $\cC$ is an extension-closed subcategory of a triangulated category $\cA$. Statements (3) and (4) follow from some well-known properties of $\cA$; see \cite[(1.4)]{H} and \cite[(2.4),(2.5)]{Kra2}. Given a morphism $u: X\to M$ in $\cC$, we obtain a commutative diagram with rows being exact triangles \vspace{-2pt}
$$\xymatrixrowsep{16pt}\xymatrix{X\ar[r]^f \ar[d]_-u & Y\ar[r]^g \ar@{.>}[d] & Z\ar@{=}[d] \ar[r]^\delta & X[1] \ar[d]^{u[1]}\\
M\ar[r]^{f'} & L \ar[r]^{g'} & Z \ar[r]^-{u\cdot \delta} & M[1]}$$ in $\cA$, where $L\in \cC$.
If $u\cdot \delta=0,$ then $f'$ is a section, and hence, $u$ factors through $f$. If $u$ factors through $f$ then, by rotating the top exact triangle to the left, we see that $u\circ (-\delta[-1])=0$, and hence, $u\cdot \eta=u[1]\circ \delta=0$. This establish Statement (1). Dually, we can prove Statement (2). The proof of the lemma is completed.

\medskip

We are ready to introduce the stable categories of $\cC$. Fix an object $M$ in $\cC$. Given a morphism $f: X\to Y$ in $\cC$, we obtain two $\Z$-linear maps
$$\Ext^1_\mathcal{C}(M, f): \Ext^1_\mathcal{C}(M, X)\to \Ext^1_\mathcal{C}(M, Y):  \delta\mapsto f \cdot  \delta \vspace{-5pt} $$
and \vspace{-2pt}
$$\Ext^1_\mathcal{C}(f, M): \Ext^1_\mathcal{C}(Y, M)\to \Ext^1_\mathcal{C}(X, M): \zeta \mapsto \zeta\cdot f.\vspace{2pt} $$
This consideration yields a covariant functor $\Ext^1_\mathcal{C}(M, -): \cC\to \Mod \Z$ and a contravariant functor
$\Ext^1_\mathcal{C}(-, M): \cC\to \Mod \Z.$

\medskip

A morphism $f: X\to Y$ in $\cC$ is called {\it injectively trivial} if $\Ext_\mathcal{C}^1(M, f)=0$ for all $M\in \cC;$
and {\it projectively trivial} if $\Ext_\mathcal{C}^1(f, M)=0$ for all $M\in \cC$;
compare \cite[Section 2]{LeZ}.
Moreover, an object $X\in \cC$ is called {\it Ext-injective} if $\id_X$ is injectively trivial, or equivalently, $\Ext_\mathcal{C}^1(M, X)= 0$ for all $M\in \cC$; and {\it Ext-projective} if $\id_X$ is projectively trivial, or equivalently, $\Ext_\mathcal{C}^1(X, N)= 0$ for all $N\in \cC$. Clearly, the injectively trivial morphisms and the projectively trivial morphisms form two ideals written as $\mathcal{P}_\mathcal{C}$ and $\mathcal{I}\hspace{.3pt}_\mathcal{C}$ in $\cC$, respectively. The following observation is important.


\begin{Lemma}\label{left-tri-ideal}

Let $\cC$ be an extension-closed subcategory of a triangulated category.

\vspace{-1pt}

\begin{enumerate}[$(1)$]

\item If $X\in \cC \cap \cC[1]$, then $\mathcal{P}_\mathcal{C}(M, X)=0$ for all $M\in \cC$.

\vspace{1pt}

\item If $X\in \cC \cap \cC[-1]$, then $\mathcal{I}_{\hspace{.5pt}\mathcal{C}}(X, N)=0$ for all $N\in \cC$.

\end{enumerate}\end{Lemma}

\noindent{\it Proof.} We shall only prove Statement (1). Let $f: M\to X$ be projectively trivial, where
$X\in \cC \cap \cC[1]$. Since $X[-1]\in \cC$, we see that $\id_X \in \Ext^1_{\mathcal{C}}(X, X[-1])$ is such that ${\rm Ext}_\cC^1(f, X[-1])(\id_X)=0,$ that is, $f= 0$. The proof of the lemma is completed.

\medskip

We are ready to define the stable categories of a tri-exact category.

\begin{Defn}

Let $\cC$ be a tri-exact category. We shall call $\overline{\hspace{-1pt}\cC}=\cC/\mathcal{I}_\mathcal{C}$ the {\it injectively stable category}, and $\underline{\cC\hspace{-1pt}}=\cC/\mathcal{P}_\mathcal{C}$ the {\it projectively stable category}, of $\cC.$

\end{Defn}

\medskip

\noindent{\sc Remark.} In view of Lemmas \ref{left-tri} and \ref{left-tri-ideal}, we see that $\underline{\cC\hspace{-1pt}}=\cC$ in case $\cC$ is a left triangulated category, and $\oC=\cC$ in case $\cC$ is a right triangulated category.

%

\medskip

We shall put \vspace{.5pt} $\oHom_\mathcal{\hspace{.8pt}C}(X, Y)= \Hom_{\oC}(X, Y)$ and $\uHom_\mathcal{\hspace{1pt}C}(X, Y)=\Hom_{\hspace{1pt}\underline{\mathcal C\hspace{-1pt}}\hspace{1pt}}(X, Y)$, for all $X, Y\in \cC$. Consider a morphism $f: X\to Y$ in $\cC$. We shall write $\bar{\hspace{.5pt}f}$ and $\underline{f\hspace{-2pt}}\hspace{2pt}$ for its images in
$\overline{\rm Hom}_\mathcal{\hspace{.8pt}C}(X, Y)$ and $\underline{\hspace{-.5pt}\rm Hom\hspace{-.5pt}}_\mathcal{\hspace{2pt}C}(X, Y)$, respectively. In this way, we may define $\bar{\hspace{.5pt}f} \cdot \zeta=f \cdot \zeta$ and $\delta \cdot \underline{f\hspace{-2pt}}\hspace{2pt}=\delta \cdot f$, for all $\zeta\in \Ext^1(M, X)$ and $\delta\in \Ext^1(M, Y)$.

\begin{Lemma}\label{Lift-stab-obj}

Let $\cC$ be a Krull-Schmidt tri-exact category.

\vspace{-2pt}

\begin{enumerate}[$(1)$]

\item If $X\in \uC\hspace{.5pt}$ is indecomposable, then there exists an indecomposable object $M$ in $\cC$
such that $\uHom_\sC(X, -)\cong \uHom_\sC(M, -)$ and $\Ext^1_\sC(X, -)\cong \Ext^1_\sC(M, -)$.

\vspace{1pt}

\item If $X\in \oC$ is indecomposable, then there exists an indecomposable object $N$ in $\cC$ such that
 $\oHom_\sC(-, X)\cong \oHom_\sC(-, N)$ and $\Ext^1_\sC(-, X)\cong \Ext^1_\sC(-, N)$.

\end{enumerate}

\end{Lemma}

\noindent{\it Proof.} We shall only prove the first statement. Let $X\in \uC$ be indecomposable. Since $\cC$ is Krull-Schmidt, $\End(X)$ is semiperfect; see \cite[(1.1)]{LNP}. Thus, $\End(X)$ has a complete orthogonal set $\{e_1, \ldots, e_n\}$ of primitive idempotents such that $e_i \End(X) e_i$ is local, for $i=1, \ldots, n$; see \cite[(27.6)]{AnF}. Since $\underline{\End\hspace{-1pt}}\hspace{1.5pt}(X)$ is local, we may assume that $\underline{\hspace{-.5pt}e\hspace{-.8pt}}\hspace{.5pt}_1=\underline{\id}\hspace{1pt}_X.$ Let $q: M\to X$ and $p: X\to M$ be morphisms such that $p \hspace{.5pt}q=\id_M$ and $q p=e_1$. Observing that $\End(M)\cong e_1 \End(X) e_1$, we see that $M$ is indecomposable in $\cC$. Given $Y\in \cC$, since $\id_X-e_1$ is projectively trivial, we obtain two $\Z$-linear isomorphisms $\Ext^1_\cC(p\hspace{.5pt}, Y): \Ext_\cC^1(M, Y)\to \Ext_\cC^1(X, Y)$ and $\uHom_\sC(\hspace{1pt}\underline{p\hspace{-1pt}}\hspace{1.5pt}, Y): \uHom_\sC(M, Y)\to \uHom_\sC(X, Y)$, which are evidently natural in $Y$.
The proof of the lemma is completed.

\medskip

Next, we shall relate exact categories to tri-exact categories. Fix an abelian category $\mf A$ and consider its derived category $D(\mf A)$. Given an object $X$ in $\mf{A}$, we obtain a stalk complex $X[n]$ whose component of degree $-n$ is $X$. Given a morphism $f: X\to Y$ in $\mf{A}$, we obtain a morphism $f[n]: X[n]\to Y[n]$ whose component of degree $-n$ is $f$. Consider the canonical embedding functor
$$\mathbb{D}: \mf{A}\to D(\mf{A}): X\mapsto X[0]; f\mapsto \tilde{f}[0], \vspace{-1pt}$$ where $\tilde{f}[0]$ is the image of $f[0]$ under $\mathbb{L}\circ \mathbb{P}: C(\mf{A})\to D(\mf{A})$; see \cite[(III.3.4.7)]{Mil}. Let $\mathscr{C}$ be an extension-closed subcategory of $\mf A$. We shall denote by $\hspace{-2pt}\hat{\hspace{2pt}\mathscr{C}}$ the full subcategory of $D(\mf{A})$ of complexes $X^\cdt$ with ${\rm H}^0(X^\cdt)\in \mathscr{C}$ and ${\rm H}^i(X^\cdt)=0$ for  $i\ne 0$.


\begin{Lemma}\label{ext-closed}

Let $\mathscr{C}$ be an extension-closed subcategory of an abelian category $\mf{A}$.

\vspace{-4pt}

\begin{enumerate}[$(1)$]

\item The category $\hspace{-2pt}\hat{\hspace{2pt}\mathscr{C}}$ is an extension-closed subcategory of $D(\mf{A})$.

\item Given $X^\cdt\in \hspace{-2pt}\hat{\hspace{2pt}\mathscr{C}}$, there exists a natural isomorphism $\theta_{\hspace{-.5pt}X^\cdt}: X^\cdt\to {\rm H}^0(X^\cdt)[0]$ in $\hspace{-2pt}\hat{\hspace{2pt}\mathscr{C}}.$

\vspace{-1pt}

\end{enumerate}

\end{Lemma}

\noindent{\it Proof.} Given $X^\cdt\in \hspace{-2pt}\hat{\hspace{2pt}\mathscr{C}}$, consider the smart truncations $\tau_{\le 0}(X^\cdt)$ and $\tau_{\ge 0}(\tau_{\le 0}(X^\cdt))$. Since ${\rm H}^i(X^\cdt)=0$ for all $i\ne 0$, the canonical injection $i_0^\sdt: \tau_{\le 0}(X^\cdt) \to X^\cdt$ and the canonical projection $p_0^\sdt: \tau_{\le 0}(X^\cdt) \to \tau_{\ge 0}(\tau_{\le 0}(X^\cdt))\vspace{1pt}$ are quasi-isomorphisms; see \cite[(III.3.4)]{Mil}. Observing that $\tau_{\ge 0}(\tau_{\le 0}(X^\cdt))={\rm H}^0(X^\cdt)[0]$, we obtain a isomorphism $\theta_{\hspace{-.5pt}X^\pdt}: X^\cdt\to {\rm H}^0(X^\cdt)[0]$ in $\hspace{-2pt}\hat{\hspace{2pt}\mathscr{C}}$, which is evidently natural in $X^\ydt$.

Let $\hspace{-3pt}\xymatrixcolsep{18pt}\xymatrix{ X^\ydt \ar[r] & Y^\ydt \ar[r] & Z^\ydt \ar[r] & X^\ydt[1]\hspace{-3pt}},$\vspace{-3pt} where  $X^\ydt, Z^\ydt\in \hspace{-2pt}\hat{\hspace{2pt}\mathscr{C}}$, be an exact triangle in $D(\mf{A})$. By the long exact sequence of cohomology, 
${\rm H}^i(Y^\cdt)=0$  \vspace{-2pt} for all $i\ne 0$ and $\hspace{-2pt}\xymatrixcolsep{18pt}\xymatrix{0\ar[r] & {\rm H}^0(X^\ydt) \ar[r] & {\rm H}^0(Y^\ydt) \ar[r] & {\rm H}^0(Z^\ydt) \ar[r] & 0}\hspace{-2pt}\vspace{-1.5pt}$ is a short exact sequence in $\mf{A}$. Since ${\rm H}^0(X^\ydt), {\rm H}^0(Z^\ydt)\in \mathscr{C}\vspace{1pt}$, we see that ${\rm H}^0(Y^\pdt) \in \mathscr{C}$. The proof of the lemma is completed.

\medskip

By Lemma \ref{ext-closed}, restricting the canonical embedding $\mathbb{D}: \mf{A}\to D(\mf{A})$ yields an equivalence $\mathbb{D}_\mathscr{\hspace{.4pt}C}: \mathscr{C} \to \hat{\hspace{1.5pt}\mathscr{C}}.$ We shall say that a functor $F: \mathscr{C}\to \Mod \Z$ is {\it essentially equivalent} to a functor $\hat{F}: \hspace{-2pt}\hat{\hspace{2pt}\mathscr{C}} \to \Mod \Z$ provided that $F\cong \hat{F}\circ \mathbb{D}\hspace{.4pt}_\mathscr{C}.$

\begin{Lemma}\label{Fun-equiv}

Let $\mathscr{C}$ be an extension-closed subcategory of an abelian category $\mf{A}$, and consider an object $X$ in $\mathscr{C}$.

\begin{enumerate}[$(1)$]

\item The functors $\Ext_\mathscr{C}^1(X, -)$ and $\Ext_\mathscr{C}^1(-, X)$ are essentially equivalent to the functors $\Ext_{\hspace{-1pt}\hat{\hspace{2pt}\mathscr{C}}}^1(X[0], -)\vspace{1.5pt}$ and $\Ext_{\hspace{-1pt}\hat{\hspace{2pt}\mathscr{C}}}^1(-, X[0]),$ respectively.

\vspace{0pt}

\item The functors $\uHom_\mathscr{\hspace{.5pt}C}(X, -)\vspace{1.5pt}$ and $\oHom_\mathscr{\hspace{.5pt}C}(\hspace{-1pt}-, X)$ are essentially equivalent to the functors $\uHom_{\hspace{-1pt}\hat{\hspace{2pt}\mathscr{C}}}(X[0], \hspace{-1pt}-\hspace{-1pt})$ and $\oHom_{\hspace{-2pt}\hat{\hspace{2pt}\mathscr{C}}}(-, \hspace{-1pt}X[0]),$ respectively.

\end{enumerate}\end{Lemma}

\noindent{\it Proof.} Given an extension $\delta\in \Ext_{\hspace{.4pt}\mathscr{C}}^1(M, N)$ represented by a short exact sequence \vspace{-3pt} $$\xymatrixcolsep{18pt}\xymatrix{0\ar[r] & N \ar[r]^f& L \ar[r]^g & M \ar[r] &0}$$ in $\mathscr{C}$, it is well-known that $D(\mf A)$ has an induced exact triangle \vspace{-3pt}
$$\xymatrix{N[0] \ar[r]^{\tilde{f}[0]} & L[0] \ar[r]^{\tilde{g}[0]} & M[0] \ar[r]^{\tilde\delta} & N[1].}
\vspace{-2pt}$$
This yields an isomorphism $$\mathbb{E}_{M, N}: \Ext_{\hspace{.4pt}\mathscr{C}}^1(M, N) \to \Ext_{\hspace{-1pt}\hat{\hspace{2pt}\mathscr{C}}}^1(M[0], N[0]): \delta\mapsto \tilde{\delta},$$ such, for all $f\in \Hom_\mathscr{\hspace{.5pt}C}(X, M)$,
$\delta\in \Ext_\mathscr{\hspace{.5pt}C}^1(Z,X)$ and $g\in \Hom_\mathscr{\hspace{.5pt}C}(N, Z)$, that $\mathbb{E}_{N,M}(g\cdot \delta \cdot f)=\tilde{f}[0]\cdot \tilde{\delta}\cdot \tilde{g}[0]$; see \cite[(IV.2.1.1)]{Mil}.
Thus, Statement (1) holds.

\vspace{.5pt}

To show Statement (2), we claim that a morphism $f: X\to Y$ in $\mathscr{C}$ is injectively trivial if and only if $\tilde{f}[0]$ is injectively trivial in $\hspace{-2pt}\hat{\mathscr{\hspace{2pt}C}}$. Assume first that \vspace{.5pt} $\tilde{f}[0]$ is injectively trivial in $\hspace{-2pt}\hat{\mathscr{\hspace{2pt}C}}$. Given any $\delta\in \Ext^1_\mathscr{\,C}(Z, X)$, we obtain $\mathbb{E}_{Z,Y}(f\cdot \delta)=\tilde{f}[0] \cdot \tilde{\delta}=0$, and hence, $f\cdot \delta=0$. That is, $f$ is injectively trivial in $\mathscr{C}$. \vspace{.3pt} Conversely, assume that $f$ is injectively trivial in $\mathscr{C}.$ \vspace{-1pt} Let $\zeta^\sdt\in \Ext_{\hspace{-2pt}\hat{\mathscr{\hspace{2pt}C}}}^1(Z^\pdt, X[0])$. Setting $Z={\rm H}^0(Z^\pdt)$, by Lemma \ref{ext-closed}(2), we have an isomorphism $\theta_{Z^\pdt}: Z^\pdt\to Z[0]$ in $\hspace{-2pt}\hat{\mathscr{\hspace{2pt}C}}$. Thus, $\zeta^\sdt \cdot \theta_{\hspace{-1pt}Z^\pdt}^{-1}=\tilde{\delta}$ for some $\delta\in \Ext_\mathscr{C}^1(Z, X)$. Observing that
$\tilde{f}[0]\cdot \tilde{\delta}=\mathbb{E}_{Z,Y}(f \cdot \delta)=\mathbb{E}_{Z,Y}(0)=0,$ we obtain
$\tilde{f}[0]\cdot \zeta^\sdt=0$. That is, $\tilde{f}[0]$ is injectively trivial in $\hspace{-2pt}\hat{\hspace{2pt}\mathscr{C}}$. This establishes our claim. As a consequence, $\mathbb{D}_\mathscr{C}: \mathscr{C}\to \hat{\hspace{2pt}\mathscr{C}}$ induces an equivalence between the injectively stable categories. In particular, $\uHom_\mathscr{\hspace{.5pt}C}(X, -)\vspace{1.5pt}$
is essentially equivalent to $\uHom_{\hspace{-1pt}\hat{\hspace{2pt}\mathscr{C}}}(X[0], \hspace{-1pt}-\hspace{-1pt})$. In a dual fashion, we may establish the second part of Statement (2). The proof of the lemma is completed.

\medskip

The following statement is needed in the next section.


\begin{Lemma} \label{Fun-mono}

Let $\mathscr{C}$ be an extension-closed subcategory of an abelian category $\mf{A}$, and let $F, G: \mathscr{C}\to \Mod \Z$ be
functors essentially equivalent to $\hat{F}, \hat{G}: \hspace{-2pt}\hat{\hspace{2pt}\mathscr{C}} \to \Mod \hspace{.3pt} \Z$ respectively. Then there exists a $(\hspace{-1.5pt}$mono, iso$)$morphism $\eta:  F\to G$ if and only if there exists a $(\hspace{-1.5pt}$mono, iso$)$morphism $\hat{\eta}: \hat{F}\to \hat{G}$.

\end{Lemma}

\noindent{\it Proof.} Let  $\zeta: F\to \hat{F} \circ \mathcal{D}_\mathscr{C}$ and $\xi: G\to \hat{G} \circ \mathcal{D}_\mathscr{C}$ be isomorphisms. Firstly, assume that $\hat{\eta}: \hat{F}\to \hat{G}$ is a morphism. Given $X\in \mathscr{C}$, we
set $\eta_X=\xi_X^{-1}\circ \hat{\eta}_{X[0]}\circ \zeta_X$, which is a monomorphism or isomorphism in case $\hat{\eta}_{X[0]}$ is a monomorphism or isomorphism, respectively. This yields a desired (mono, iso)morphism $\eta:  F\to G$.

Conversely, assume that $\eta: F \to G$ is a morphism. Given $X^\cdt\in \hspace{-2pt}\hat{\hspace{2.5pt}\mathscr{C}}$, by Lemma \ref{ext-closed}(2), there exists a natural isomorphism $\theta_{\hspace{-.5pt}X^\pdt}: X^\cdt\to X[0]$, where $X={\rm H}^0(X^\cdt).$ We define $\hat{\eta}_{X^\cdt}$ to be the composite of the following morphisms \vspace{-3pt}
$$\xymatrix{\hat{F}(X^\cdt)\ar[r]^{\hat{F}(\theta_{X^\pdt})} & \hat{F}(X[0])\ar[r]^{\zeta_X^{-1}} &
F(X)\ar[r]^{\eta_X}& G(X) \ar[r]^{\xi_X}& \hat{G}(X[0]) \ar[r]^{\hat{G}(\theta_{X^\cdt}^{-1})} & \hat{G}(X^\cdt),}$$
which is a monomorphism or isomorphism if $\eta_X$ is a monomorphism or isomorphism respectively. This yields a desired (mono, iso)morphism $\hat{\eta}: \hat{F}\to \hat{G}$. The proof of the lemma is completed.


\section{Almost split sequences}

\medskip

\noindent The objective of this section is to study the existence of an individual almost split sequence in a tri-exact category. Using
similar but more general techniques, we shall unify and extend the results under various classical settings; see \cite{ABM, Kra2, LNP, LeZ, RvdB}. In parti\-cular, Auslander's existence theorem for an almost split sequence in the category of all modules over a ring; see \cite{AUS} and Krause's existence theorem for an almost triangle in a triangulated category fit well into our setting; see \cite{Kra2}.

\medskip

Throughout this section, $\cC$ stands for a tri-exact category, say an extension-closed subcategory of a triangulated category $\cA$. The following notion plays a fundamental role in our investigation.


\begin{Defn}\label{art}

Let $\cC$ be a tri-exact category. A tri-exact sequence $$\xymatrix{X \ar[r]^f & Y \ar[r]^g &Z}$$ defined by an extension $\delta\in \Ext_\mathcal{C}^1(Z,  X)$ is called {\it almost split} if $f$ is minimal left almost split and $g$ is minimal right almost split. In this case, $\delta$ is called {\it almost-zero}.

\end{Defn}

\medskip

\noindent{\sc Remark.} An almost split sequence with a non-zero middle term is an Auslander-Reiten sequence as defined in \cite[(1.3)]{Liu1}. The converse is probably not true.

\medskip

Modifying slightly the proof of the proposition stated in \cite[(3.5)]{Hap2}, we obtain the uniqueness of an almost split sequence in a tri-exact category as follows.

\begin{Prop}\label{ARS-uni}

Let $\cC$ be a tri-exact category. If $\hspace{-3pt}\xymatrixcolsep{18pt}\xymatrix{X\ar[r]&Y\ar[r] &Z}\hspace{-3pt}$ is an almost split sequence in $\cC$, then it is unique up to isomorphism for $X$ and for $Z$.

\end{Prop}

\medskip

The following statement says in particular that the study of almost split sequences under various classical settings can be unified under our tri-exact setting.


\begin{Prop}\label{ART-ARS}

Let $\mathscr{C}$ be an extension-closed subcategory \vspace{-1pt} of an abelian category $\mf{A}$. Then every almost split sequence $\xymatrixcolsep{18pt}\xymatrix{0\ar[r] & X \ar[r] & Y\ar[r]  & Z \ar[r] & 0}\vspace{-2pt}$ in $\mathscr{C}$ induces an almost split sequence $\xymatrixcolsep{18pt}\xymatrix{\hspace{-3pt} X[0] \ar[r] & Y[0] \ar[r]  & Z[0]\hspace{-3pt}}\vspace{-3pt}$ in $\hspace{-2pt}\hat{\mathscr{\hspace{2pt}C}};$ and every almost split sequence in $\hspace{-2pt}\hat{\mathscr{\hspace{2pt}C}}$ can be obtained in this way.

\end{Prop}

\noindent{\it Proof.} Since $\mathbb{D}_\mathscr{C}: \mathscr{C}\to \hspace{-2pt}\hat{\hspace{2pt}\mathscr{C}}$ is an equivalence, the first part of the proposition follows immediately. Assume that $\hspace{-2pt}\hat{\mathscr{\hspace{2pt}C}}$ has an almost split sequence \vspace{-2pt} $$\xymatrix{\hspace{-3pt} X^\ydt\ar[r]^{f^\sdt} & Y^\ydt\ar[r]^{g^\ydt} & Z^\ydt\hspace{-3pt}}$$ defined by an extension $\eta^\ydt\in {\rm Ext}_{\hspace{-2pt}\hat{\mathscr{\hspace{2pt}C}}}^1(Z^\ydt, X^\ydt[1])$. Applying the long exact sequence of cohomology, we obtain a short exact sequence \vspace{-2pt}
$$\delta: \quad \xymatrix{0\ar[r] & {\rm H}^0(X^\ydt)\ar[r]^{f} & {\rm H}^0(Y^\ydt)\ar[r]^{g} &
{\rm H}^0(Z^\ydt)\ar[r] & 0} \vspace{-1pt}$$ in $\mathscr{C},$ where $f={\rm H}^0(f^\ydt)$ and $g={\rm H}^0(g^\ydt)$. In view of Lemma \ref{ext-closed}(2), there exists a commutative diagram with vertical isomorphisms \vspace{-2pt}
$$\xymatrixrowsep{16pt}\xymatrix{X^\ydt\ar[r]^{f^\sdt} \ar[d]_-{\theta_{\hspace{-.5pt}X^\ydt}} & Y^\ydt \ar[r]^{g^\ydt} \ar[d]^-{\theta_{\hspace{-.5pt} Y^\ydt}} & Z^\ydt \ar@{.>}[d]^-{\zeta^\sdt} \ar[r]^{\eta^\ydt} & X^\cdt[1] \ar[d]^-{\theta_{\hspace{-1pt}X^\ydt}[1]}\\
{\rm H}^0(X^\ydt)[0] \ar[r]^{\tilde{f}[0]} & {\rm H}^0(Y^\ydt)[0]\ar[r]^{\tilde{g}[0]} & {\rm H}^0(Z^\ydt)[0]
\ar[r]^{\tilde{\delta}} & {\rm H}^0(X^\ydt)[1]
}$$
in $D(\mf{A})$, where the rows are exact triangles. In particular, $\tilde{f}[0]$ is minimal left almost split and $\tilde{g}[0]$ is minimal right almost split in $\hat{\hspace{2pt}\mathscr{C}}$. Since $\mathbb{D}_\mathscr{C}: \mathscr{C}\to \hspace{-2pt}\hat{\hspace{2pt}\mathscr{C}}$ is an equivalence, $\delta$ is almost split in $\mathscr{C}$. The proof of the proposition is completed.

\medskip

The following characterization of an almost split sequence in a tri-exact category is adapted from those in the classical settings; see \cite{AuR2, Kra2}.


\begin{Theo}\label{ART-equivalence}

Let $\cC$ be a tri-exact category. If $\xymatrix{X\ar[r]^f &Y\ar[r]^g &Z}$ is a tri-exact sequence in $\cC$, then the following statements are equivalent. \vspace{-3pt}

\begin{enumerate}[$(1)$]

\item The sequence is an almost split sequence in $\cC$.

\item The morphism $f$ is left almost split and $g$ is right almost split.

\item The morphism $f$ is left almost split and $Z$ is strongly indecomposable.

\item The morphism $g$ is right almost split and $X$ is strongly indecomposable.

\item The morphism $f$ is minimal left almost split or $g$ is minimal right almost split.

\end{enumerate}

\end{Theo}

\noindent{\it Proof.} Let $\eta\in \Ext_\cC^1(Z, X)$, defining the tri-exact sequence stated in the theorem.
If Statement (2) holds, then $X$ and $Z$ are strongly indecomposable; see \cite[(2.3)]{AuR2}, and by Lemma \ref{section-retraction}(4), $g$ is right minimal and $f$ is left minimal, that is, Statement (1) holds. Moreover,  by Lemma \ref{section-retraction}(4),
either of Statements (3) and (4) implies Statement (5). Assume now that Statement (5) holds. To prove Statement (2), we assume that $\cC$ is an extension-closed subcategory of a triangulated category $\cA$. If $g$ is minimal right almost split, using the same argument given in \cite[(2.6)]{Kra2}, we may show that $f$ is left almost split. If $f$ is minimal left almost split, one can dually show that $g$ is right almost split.
The proof of the theorem is completed.

\medskip

The rest of this section is devoted to the study of the existence of an almost split sequence in a tri-exact category. We start with some properties of almost-zero extensions; compare \cite[(2.2)]{JOR}, \cite[(3.1)]{LeZ} and \cite[Page 306]{RvdB}.


\begin{Lemma}\label{art-1}

Let $\cC$ be a tri-exact category with $\delta\in \Ext^1_\mathcal{C}(Z, X)$ being almost-zero.

\vspace{-1pt}

\begin{enumerate}[$(1)$]

\item A factorization $\delta = \eta \cdot  \underline{f\hspace{-2pt}}\hspace{2pt}$ exists whenever
a non-zero extension $\eta\in \Ext^1_\mathcal{C}(L, X)$ or a non-zero morphism $\underline{f\hspace{-2pt}}\hspace{2pt} \in \uHom_{\hspace{.6pt}\mathcal{C}}(Z, L)$ is given.

\vspace{0pt}

\item A factorization $\delta=\bar{\hspace{-.6pt}g} \cdot \zeta$ exists whenever a non-zero extension $\zeta \in \Ext^1_{\mathcal{C}}(Z, L)$ or a non-zero morphism $\bar{\hspace{-.6pt}g}\in \oHom_{\hspace{.5pt}\mathcal{C}}(L, X)$ is given.

\end{enumerate}

\end{Lemma}

\noindent{\it Proof.} 
Given a non-zero extension $\eta\in \Ext^1_\mathcal{C}(L, X)$, using the same proof of the first statement of the sublemma stated in \cite[Page 306]{RvdB}, \vspace{-2pt} we obtain $\delta = \eta \cdot  \underline{f\hspace{-2pt}}\hspace{2pt}$ for some morphism $\underline{f\hspace{-2pt}}\hspace{2pt} \in \uHom_{\hspace{.6pt}\mathcal{C}}(Z, L)$. Given a non-zero extension $\zeta \in \Ext^1_{\mathcal{C}}(Z, L)$, by a dual argument, we can show that $\delta=\bar{\hspace{-.6pt}g} \cdot \zeta$ for some $\bar{\hspace{-.6pt}g}\in \oHom_{\hspace{.5pt}\mathcal{C}}(L, X)$.

Now, let $f: Z\to L$ be a non projectively trivial morphism in $\cC$. Then, there exists some $\zeta\in
\Ext_\cC(L, M)$ such that $0\ne \zeta\cdot f\in \Ext_\cC(Z, M)$. By the first part of Statement (2), there exists some $g: M\to X$ in $\cC$ such that $\delta= g\cdot (\zeta \cdot f)=(g \cdot \zeta) \cdot f$. This establishes the second part of Statement (1). Dually, one can verify the second part of Statement (2). The proof of the lemma is completed.

\medskip

Given $X, Z\in \cC$, by Lemma \ref{Ext-bimod}, $\Ext_\mathcal{C}^1(Z, X)$ is an $\overline{\End}(X)$-$\underline{\End\hspace{-1pt}}\hspace{1pt}(Z)$-bimodule so that $\overline{\hspace{-2pt}f} \cdot \delta \cdot \underline{g\hspace{-1pt}} = f\cdot \delta \cdot g$, for $f\in \End(X), \delta \in \Ext_\mathcal{C}^1(Z, X)$ and $g\in \End(Z).$ Given
$M\in \cC$ strongly indecomposable, we always write $S_M=\End(M)/{\rm rad}(\End(M))$, which is a simple left $\overline{\End}(M)$-module if $M$ is not Ext-injective; and a simple right $\underline{\End\hspace{-.8pt}}\,(M)$-module if $M$ is not Ext-projective.

\begin{Theo}\label{AZE}

Let $\mathcal{C}$ be a tri-exact category with $X, Z\in \cC$ strongly indecomposable. Consider non-zero ring homomorphisms $\Ga\to \overline{\End}(X)$ and \vspace{0pt} $\Sa\to \underline{\End}(Z)$. Let $_{\it\Gamma}\hspace{-.4pt}I$ be
an injective cogenerator of the left $\Ga$-module $S_X\vspace{1pt}$ and $I_{\hspace{-.5pt}\it\Sigma}$ an injective cogenerator of the right $\Sa$-module \vspace{.5pt} $S_Z$. If $\delta\in \Ext^1_\mathcal{C}(Z, X)$ is non-zero, then $\delta$ being almost zero is equivalent to each of the following statements.

\vspace{0pt}

\begin{enumerate}[$(1)$]


\item There exists a monomorphism $\Psi: \Ext_\mathcal{C}^1(Z,-)\to \Hom{_{\it\Gamma}}(\hspace{.8pt}\oHom_\mathcal{\hspace{1pt}C}(-, X), {}_{\it\Gamma}\hspace{-.4pt}I)\vspace{.5pt}$ such that $\Psi_X(\delta)$ lies in the socle of the left $\overline{\End}(X)$-module $\Hom{_{\it\Gamma}}(\hspace{.8pt}\overline{\rm End}(X), {}_{\it\Gamma}\hspace{-.4pt}I)$.

\vspace{0.5pt}

\item There exists a monomorphism $\Phi: \Ext^1_\mathcal{C}(-, X)\to \Hom_{\it\Sigma}(\uHom_\mathcal{\hspace{1.5pt}C}(Z,-), I_{\hspace{-.6pt}\it\Sigma})\vspace{.5pt}$
    such that $\Phi_Z(\delta)$ lies in the socle of the right $\underline{\End\hspace{-.5pt}}\hspace{.4pt}(\hspace{-.5pt}Z)$-module $\Hom_{\it\Sigma}(\hspace{.5pt}\underline{\rm End\hspace{-1pt}}(Z), I_{\hspace{-.6pt}\it\Sigma})$.

\end{enumerate}

\end{Theo}

\noindent{\it Proof.} Let $\delta\in \Ext^1_{\cC}(Z, X)$ be non-zero. We shall only consider Statement (1). Let $I$ be an injective cogenerator of the left $\Ga$-module $S_X$. Assume first that $\delta$ is almost-zero. Consider the $\overline{\End}(X)$-submodule $S$ of $\Ext^1_{\cC}(Z, X)$ generated by $\delta$, which is simple by Lemma \ref{art-1}(2). Since $\overline{\End}(X)$ is local, $S\cong S_X$ as left $\overline{\End}(X)$-modules, and consequently, $S\cong S_X$ as $\Ga$-modules.
Thus, we may find a $\Ga$-linear map $\psi: \Ext^1_\mathcal{C}(Z, X) \to I$ with $\psi(\delta)\neq 0.$
Given $L\in \mathcal{C}$, by Lemma \ref{art-1}(2), we obtain a non-degenerate $\Z$-bilinear form
$$<\hspace{-3pt}-\,,-\hspace{-3pt}>_L \hspace{1pt}: \hspace{2pt} \overline{\rm Hom}_\mathcal{\hspace{.5pt}C}(L,X) \times \Ext^1_\mathcal{C}(Z, L)\to I:(\bar{g}, \zeta)\mapsto \psi(\bar{g}\cdot \zeta).$$
Since $\psi$ is left $\Ga$-linear, by Lemma \ref{Ext-bimod}, we obtain a $\Z$-linear monomorphism $$\Psi_L: \Ext^1_\mathcal{C}(Z,L) \to \Hom_{\it\Gamma}(\hspace{.7pt}\overline{\rm Hom}_\mathcal{\hspace{.5pt}C}(L, X), I): \zeta\mapsto <\hspace{-3pt}-\,,\zeta\hspace{-1pt}>_L,$$
which is natural in $L$. This yields a monomorphism $\Psi$ as stated in Statement (1).
Since $\Psi_X$ is a left $\overline{\End}(X)$-linear monomorphism, $\Psi_X(S)$ is a simple submodule of the left $\overline{\End}(X)$-module $\Hom_{\it\Gamma}(\hspace{.5pt}\overline{\rm End}(X), I).$ This establishes Statement (1).

\vspace{2.5pt}

Conversely, assume that $\Psi: \Ext^1_\mathcal{C}(Z,-) \to \Hom_{\it\Gamma}(\,\overline{\rm Hom}_\mathcal{\hspace{.5pt}C}(-, X), I)\vspace{1.5pt}$ is a monomorphism such that $\Psi_X(\delta)$ is in the left $\overline{\rm End}(X)$-socle of $\Hom_{\it\Gamma}(\hspace{.5pt}\overline{\rm End}(X), I).$ Then, $\Psi_X(\delta)$ vanishes on
${\rm rad}(\overline{\End}(X)).$ Consider the non-splitting tri-exact sequence \vspace{-4pt}
$$(*) \qquad
\xymatrix{X\ar[r]^f &Y\ar[r]^g &Z} \vspace{-4pt}$$
in $\cC$ defined by $\delta$. Let $u: X \to L$ be a non-section morphism in $\mathcal{C}$. In view of the commutative diagram  \vspace{-4pt}
$$\xymatrixcolsep{30pt}\xymatrixrowsep{18pt}\xymatrix{\Ext^1_\mathcal{C}(Z, X) \ar[r]^-{\Psi_X} \ar[d]_{\Ext^1_\mathcal{C}(Z, \bar{u})} & \Hom_{\it\Gamma}(\,\overline{\End}(X), I)\ar[d]^-{\Hom_{\it\Gamma}(\oHom_{\hspace{.5pt}\mathcal{C}}(\hspace{-1pt}\bar{\hspace{1pt}u}, X), I)} \\
\Ext^1_\mathcal{C}(Z, L) \ar[r]^-{\Psi_L}
& \Hom_{\it\Gamma}(\hspace{.5pt}\oHom_\mathcal{\hspace{.8pt}C}(L, X), I),}$$ we obtain $\Psi_L(\bar{u}\cdot \delta)=\Psi_X(\delta) \circ \oHom_\mathcal{\hspace{.6pt}C}(\bar{u}, X).$ For any morphism $v: L\to X$ in $\cC$, since $v u \in {\rm rad}(\End(X)),$ we see that
$\Psi_L(\bar{u}\cdot \delta)(\bar{v})=\Psi_X(\delta)(\bar{v}\hspace{.4pt}\bar{u})= 0.$
Thus, $\Psi_L(\bar{u} \cdot \delta)=0$, and hence, $\bar{u}\cdot \delta=0$. By Lemma \ref{section-retraction}(1), $u$ factors through $f$. That is, $f$ is left almost split. Since $Z$ is strongly indecomposable, by Theorem \ref{ART-equivalence}(3), the tri-exact sequence $(*)$ is almost split. The proof of the theorem is completed.

\medskip

\noindent{\sc Remark.} In view of Lemmas \ref{Fun-equiv} and \ref{Fun-mono} and Proposition \ref{ART-ARS}, it is easy to see that Theorem \ref{AZE} covers the result stated in \cite[(2.2)]{LNP}.

\medskip

We are ready to obtain our main existence theorem of an almost splits sequence.


\begin{Theo}\label{existence}

Let $\mathcal{C}$ be a tri-exact category with $X, Z\in \cC$ strongly indecomposable. Consider non-zero ring homomorphisms $\Ga\to \overline{\End}(X)$ and $\Sa\to \underline{\End}(Z)$. Let $_{\it\Gamma}\hspace{-.0pt}I$ be an injective cogenerator of the left $\Ga$-module $S_X$ and $I_{\hspace{-.5pt}\it\Sigma}$ an injective cogenerator of the right $\Sa$-module $S_Z$. The following statements are equivalent.

\vspace{-2pt}

\begin{enumerate}[$(1)$]

\item There exists an almost split sequence $\hspace{-2pt}\xymatrixcolsep{18pt}\xymatrix{X\ar[r]&Y\ar[r] &Z}\hspace{-2pt}$ in $\mathcal{C}$.

\vspace{0.5pt}

\item There exists a monomorphism
$\Psi: \Ext_\mathcal{C}^1(Z,-)\to \Hom{_{\it\Gamma}}(\hspace{.8pt}\oHom_\mathcal{\hspace{1pt}C}(-, X), {}_{\it\Gamma}\hspace{-.4pt}I)$ such that the left $\overline{\End}(X)$-linear map $\Psi_X$ is socle essential.

\vspace{0.5pt}

\item There exists a monomorphism $\Phi: \Ext^1_\mathcal{C}(-, X)\to \Hom_{\it\Sigma}(\uHom_\mathcal{\hspace{1.5pt}C}(Z,-), I_{\hspace{-.6pt}\it\Sigma})$ such that the right $\underline{\End\hspace{-.5pt}}\hspace{.4pt}(\hspace{-.5pt}Z)$-linear map $\Phi_Z$ is socle essential.

\end{enumerate}

\end{Theo}

\noindent{\it Proof.} By Theorem \ref{AZE}, it suffices to prove that Statement (2) implies Statement (1).  Let $I$ be an injective co-generator of the left $\Ga$-module $S_X$ with a monomorphism $\Psi: \Ext_\mathcal{C}^1(Z,-)\to \Hom{_{\it\Gamma}}(\hspace{.8pt}\oHom_\mathcal{\hspace{1pt}C}(-, X), I)$ such that the left $\overline{\End}(X)$-linear map $\Psi_X: \Ext_\mathcal{C}^1(Z,X)\to \Hom{_{\it\Gamma}}(\hspace{.8pt}\overline{\End\hspace{-.5pt}}(X), I)$ is socle essential. Consider the canonical projection $p: \overline{\rm End}(X) \to S_X$ and fix a non-zero $\Ga$-linear map $q: S_X\to I$. Then, $qp$ is a non-zero element in $\Hom_{\it\Gamma}(\hspace{.8pt}\overline{\rm End}(X), I)$ annihilated by ${\rm rad}(\overline{\rm End}(X))$. Since $\overline{\rm End}(X)$ is local, $qp$ belongs to the left $\overline{\rm End}(X)$-socle of $\Hom_{\it\Gamma}(\hspace{.5pt}\overline{\rm End}(X), I)$. Since $\Psi_X$ is socle essential, there exists $\delta\in \Ext_\cC^1(Z, X)$ such that $\Psi_X(\delta)$ lies in the left $\overline{\End}(X)$-socle of $\Hom{_{\it\Gamma}}(\hspace{.8pt}\overline{\rm End}(X), {}_{\it\Gamma}\hspace{-.4pt}I)$. By Theorem \ref{AZE}(1), $\delta$ defines an almost split sequence as stated in Statement (1). The proof of the theorem is completed.

\medskip

\noindent{\sc Remark.} Observe, for any ring $\Sa$, that a $\Sa$-linear monomorphism $f: M\to N$ is socle essential if $M$ has a nonzero socle or $N$ has an essential socle. Thus, we see from Lemmas \ref{Fun-equiv} and \ref{Fun-mono} and Proposition \ref{ART-ARS} that Theorem \ref{existence} includes the results stated in \cite[(4.1)]{LeZ} and  \cite[(2.3)]{LNP}.

\medskip

By abuse of terminology, a functor $F$ is called a {\it subfunctor} of another functor $G$ if there exists a monomorphism $F\to G$. We shall drop the additional hypotheses on $\Psi_Z$ and $\Phi_X$ stated in Theorem \ref{existence} in some special cases as below.


\begin{Theo}\label{existence-1}

Let $\mathcal{C}\vspace{1pt}$ be a tri-exact category with $X, Z\in \cC$ strongly indecomposable. Consider non-zero surjective ring homomorphisms $\Ga\to \overline{\End}(X)$ and $\Sa\to \underline{\End}(Z)\vspace{.5pt}$. Let $_{\it\Gamma}I$ be an injective envelope
of the left $\Ga$-module $S_X$ \vspace{.5pt} and $I_{\hspace{-.5pt}\it\Sigma}$ an injective envelope of the right $\Sa$-module $S_Z$. The following statements are equivalent.

\vspace{-2pt}

\begin{enumerate}[$(1)$]

\item There exists an almost split sequence $\hspace{-2pt}\xymatrixcolsep{20pt}\xymatrix{X\ar[r]&Y\ar[r] &Z}\hspace{-2pt}$ in $\mathcal{C}$.

\vspace{0pt}

\item $\Ext_\mathcal{C}^1(Z,-)$ is a non-zero subfunctor $\Hom{_{\it\Gamma}}(\hspace{.8pt}\oHom_\mathcal{\hspace{1pt}C}(-, X), {}_{\it\Gamma}\hspace{-.4pt}I).$

\vspace{1.3pt}

\item $\Ext^1_\mathcal{C}(-, X)$ is a non-zero subfunctor of $\Hom_{\it\Sigma}(\uHom_\mathcal{\hspace{1.5pt}C}(Z,-), I_{\hspace{-.6pt}\it\Sigma}).$

\end{enumerate}

\end{Theo}

\noindent{\it Proof.} By Theorem \ref{existence}, it suffices to prove that Statement (2) implies Statement (1). Let $\Psi: \Ext_\mathcal{C}^1(Z,-)\to \Hom{_{\it\Gamma}}(\hspace{.8pt}\oHom_\mathcal{\hspace{1pt}C}(-, X), I)\vspace{-1pt}$ be a non-zero monomorphism, where $I={}_{\it\Gamma}\hspace{-.4pt}I$. In parti\-cular, there exists some $L\in \cC$ such that $\Ext^1_\mathcal{C}(Z,L)\ne 0$. Thus, $\Psi_L(\delta)(\bar{f})\ne 0$ for some $\delta\in \Ext^1_\mathcal{C}(Z, L)$
and $f\in {\rm Hom}_\mathcal{\hspace{.5pt}C}(L, X)$. Since \vspace{-3pt}
$$\xymatrixcolsep{30pt}\xymatrixrowsep{16pt}\xymatrix{
\Ext^1_\mathcal{C}(Z, L) \ar[r]^-{\Psi_L} \ar[d]_{\Ext^1_\mathcal{C}(Z, \bar{f})} & \Hom_{\it\Gamma}(\,\overline{\rm Hom}_\mathcal{\hspace{.8pt}C}(L, X), I) \ar[d]^{\Hom_{\it\Gamma}(\oHom_{\mathcal{C}}(\bar{f}, X), I)}  \\
\Ext^1_\mathcal{C}(Z, X) \ar[r]^-{\Psi_{\hspace{-.6pt}X}} & \Hom_{\it\Gamma}(\hspace{.5pt}\overline{\End}(X), I)
\vspace{-2pt}}$$ is a commutative diagram, we obtain $\Psi_{\hspace{-.6pt}X}(\bar{f}\cdot \delta)(\overline{\id}_X)
=\Psi_L(\delta)(\bar{f})\ne 0.$\vspace{1pt}
Thus, $\Psi_{\hspace{-.6pt}X}$ is a non-zero left $\overline{\End}(X)$-linear monomorphism, which is also left $\Ga$-linear.

Since the ring homomorphism $\rho: \Ga\to \overline{\End}(X)$ is surjective, \vspace{.8pt} $S_X$ is a simple left $\Ga$-module and
$\rho^*=\Hom_{\it\Gamma}(\rho, I): \Hom_{\it\Gamma}(\,\overline{\End}(X), I)\to \Hom_{\it\Gamma}(\Ga, I)$ is a left $\Ga$-linear monomorphism. Let $\nu: \Hom_{\it\Gamma}(\Ga, I)\to I$ be the canonical $\Ga$-linear isomorphism. Then,
$\psi=\nu\circ \rho^*\circ \Psi_{\hspace{-.6pt}X}: \Ext^1_\mathcal{C}(Z, X) \to I$ is a non-zero $\Ga$-linear monomorphism. Since $S_X$ is
the essential $\Ga$-socle of $I$, it is contained in ${\rm Im}(\psi)$. In particular, $\Ext^1_\mathcal{C}(Z, X)$ has a simple left $\Ga$-submodule $S$ such that $\psi(S)=S_X$. Since $\rho$ is surjective, $S$ is a simple left $\overline{\End}(X)$-submodule of $\Ext^1_\mathcal{C}(Z, X),$ 
and hence, the left $\overline{\End}(X)$-linear monomorphism $\Psi_X$ is socle essential. By Theorem \ref{existence}(2), $\cC$ has a desired almost split sequence.  The proof of the theorem is completed.

\medskip

\noindent{\sc Remark.} (1) Let $Z\in \Mod \La$ be finitely presented, strongly indecomposable and not projective, where $\La$ is a ring. Since $\underline{\End\hspace{-.6pt}}_{\hspace{.6pt}\it\Lambda^{\rm op}}({\rm Tr}Z)^{\rm op} \cong \underline{\End\hspace{-.6pt}}_{\hspace{.6pt}\it\Lambda}(Z);\vspace{1pt}$ see \cite[Section I.3]{AUS}, we have
a surjective ring homomorphism from $\Sa=\End_{\it\Lambda^{\rm op}}({\rm Tr}Z)^{\rm op}$ onto $\underline{\End\hspace{-.6pt}}_{\hspace{.4pt}\it\Lambda \hspace{-.5pt}}(Z)$.
Let $I$ be the injective envelope of the right $\Sa$-module $S_Z$. Then, $X=\Hom_{\it\Sigma}({\rm Tr}Z, \, I)$ is a strongly indecomposable module in $\Mod\La$ such that $\Ext^1_{\it\Lambda}(-, X) \cong \Hom_{\it\Sigma}(\underline{\Hom\hspace{-.6pt}}_{\hspace{.6pt}\it\Lambda}(Z, -), I)$; see \cite[(I.11.3), (I.3.4)]{AUS}. By Theorem \ref{existence-1}, there exists an almost split sequence
$\xymatrixcolsep{16pt}\xymatrix{\hspace{-4pt} 0\ar[r] & X \ar[r] & Y \ar[r] & Z\ar[r] & 0}$ in $\Mod\La$. This is Auslander's theorem stated in \cite[(II.5.1)]{AUS}.

\vspace{1pt}

(2) Let $\mf{A}$ be a locally finitely presented Grothendieck ableian category with a finitely presented strongly indecomposable object $Z$. Consider the canonical projection $\Sa={\rm End}(Z) \to \underline{\End}(Z).$ Given any injective module $I\in \Mod \Sa^{\rm op}$, Krause  obtained a monomorphism $\Ext_\mf{\hspace{.4pt}A}^1(-, \tau_{_I}(Z)) \to \Hom_{\it\Sigma}(\uHom_\mf{\hspace{.6pt}A}(Z, -), I),$ for some $\tau_I(Z)\in \mf{A};$ see \cite[(1.2)]{Kra3}. However, it is not known whether or not $\tau_{_I}(Z)$ is strongly indecomposable even if $I$ is the injective envelope of $S_Z$. Thus, we cannot apply Theorem \ref{existence-1} to obtain an almost split sequence in $\mf{A}.$

\medskip

We can weaken the assumption that both $X$ and $Z$ are strongly indecomposable as stated in Theorem \ref{existence-1} in some special cases as below.


\begin{Theo}\label{existence-3}

Let $\mathcal{C}$ be an extension-closed subcategory of a triangulated category, and let $Z\in \cC$ be strongly indecomposable with a non-zero monomorphism \vspace{-2pt}
$$\Phi: \Hom_\mathcal{\hspace{.4pt}C}(-, X)  \to  \Hom_{\hspace{.4pt}\underline{\End}(Z)}(\hspace{.4pt}\uHom_\mathcal{\hspace{.4pt}C}(Z,-), I),\vspace{-2pt}$$ where $X\in \cC\cap \cC[1]$ and $I$ is an injective envelope of the right $\underline{\End}(Z)$-module $S_Z$. \vspace{-2pt} If $\hspace{.8pt}\Phi_{\hspace{-1.1pt}X}$ is bijective, then  $\cC$ has an almost split sequence $\hspace{-2pt}\xymatrixcolsep{18pt}\xymatrix{X[-1] \ar[r] & Y \ar[r] & Z.}$

\end{Theo}

\noindent{\it Proof.} Assume that $\hspace{.8pt}\Phi_{\hspace{-1.1pt}X}$ is bijective. Observing that $X[-1]\in \cC$ is such that $\Ext^1_\mathcal{C}(-, X[-1])=\Hom_\mathcal{\hspace{.4pt}C}(-, X)$. By Theorem \ref{existence-1}(3), it amounts to show that $\End(X)$ is local. Put $\Sa=\underline{\rm End\hspace{-.6pt}}(Z)$. By Lemma \ref{left-tri-ideal}(1), $\Hom_\mathcal{\hspace{.5pt}C}(Z, X)=\uHom_\mathcal{\hspace{.5pt}C}(Z, X)$, which is a right $\Sa$-module. Now, $\Phi_{\hspace{-.6pt}X}: \End(X)\to \Hom_{\hspace{.4pt}\it\Sigma}(\hspace{.4pt}\Hom_\mathcal{\hspace{.5pt}C}(Z, X), I)$ is a right $\End(X)$-linear isomorphism, while $\Phi_{\hspace{-.6pt}Z}: \Hom_\mathcal{\hspace{.5pt}C}(Z, X)\to \Hom_{\it\Sigma}(\Sa, I)$ is a right $\Sa$-linear monomorphism. Considering the canonical $\Sa$-linear isomorphism $\rho: \Hom_{\it\Sigma}(\Sa, I)\to I,$ we obtain a commutative diagram of surjective $\Z$-linear maps
$$\xymatrixcolsep{40pt}\xymatrixrowsep{16pt}\xymatrix{
\End_{\it\Sigma}(I) \ar[r]^-{\Hom_{\it\Sigma}(\rho, \,I)} \ar[d]_\theta & \Hom_{\it\Sigma}(\Hom_{\it\Sigma}(\Sa, I), I)\ar[d]^{\Hom_{\it\Sigma}(\Phi\hspace{-.6pt}_Z, \,I)}\\
\End(X) \ar[r]^-{\Phi\hspace{-.8pt}_X}  & \Hom_{\it\Sigma}(\Hom_\mathcal{\hspace{.4pt}C}(Z, X),I).}$$
Since $\End_{\it\Sigma}(I)$ is local; see \cite[(25.4)]{AnF}, it suffices to show that $\theta$ is a ring homomorphism. Fix arbitrarily a morphism $u: Z\to X$ in $\cC$. Given $f\in \End_{\it\Sigma}(I)$, in view of the above commutative diagram, we see that $\Phi_{\hspace{-1pt}X}(\theta(f))=f\circ \Phi_{\hspace{-1pt}Z}\circ \rho,$ and consequently, we obtain an equation
$$(1) \hspace{50pt} \Phi_{\hspace{-1pt}X}(\theta(f))(u)=f(\Phi_{\hspace{-1pt}Z}(u)(\id_{\it\Sigma})). \hspace{141pt} $$
On the other hand, considering the commutative diagram
$$\xymatrixrowsep{16pt}\xymatrix{
\End(X) \ar[d]_-{\Hom_\mathcal{C}(u, \,X)} \ar[r]^-{\Phi\hspace{-.8pt}_X} & \Hom_{\it\Sigma}(\Hom_\mathcal{C}(Z, X), I) \ar[d]^{\Hom_{\hspace{-1pt}\it\Sigma}(\Hom_{\mathcal{C}}(Z, u), I)}\\
 \Hom_{\mathcal{C}}(Z, X) \ar[r]^-{\Phi\hspace{-.6pt}_Z}  & \Hom_{\it\Sigma}(\Sa, I),
}$$
we obtain an equation
$$(2) \hspace{45pt} \Phi_{\hspace{-1pt}Z}(u)(\id_{\it\Sigma})=\Phi_X(\id_X)(u). \hspace{165pt}$$

Let $f_i\in \End_{\it\Gamma}(I)$, and write $g_i=\theta(f_i)\in \End(X)$, for $i=1, 2$. We deduce from the equation (1) that
$$\Phi_{\hspace{-1pt}X} (\theta(f_1f_2)) (u)= (f_1f_2) (\Phi_{\hspace{-1pt}Z}(u)(\id_Z))= f_1 \left[ f_2(\Phi_{\hspace{-1pt}Z}(u)(\id_Z)) \right] = f_1 \left(\Phi_{\hspace{-1pt}X}(g_2)(u) \right).$$
Since $\Phi_{\hspace{-.8pt}X}$ is right $\End(X)$-linear, combining the equations (1) and (2) yields
$$\Phi_{\hspace{-1pt}X}(g_1g_2)(u) \hspace{-1pt} = \hspace{-1pt} \Phi_{\hspace{-1pt}X}(g_1)(g_2u) \hspace{-1pt} = \hspace{-2pt} f_1 (\Phi_{\hspace{-1pt}Z}(g_2u)(\id_{\it\Sigma})) \hspace{-1pt} = \hspace{-2pt} f_1(\Phi_{\hspace{-1pt}X}(\id_X)(g_2u)) \hspace{-1pt} = \hspace{-2pt} f_1(\Phi_{\hspace{-1pt}X}(g_2)(u)).$$
Thus, $\Phi_X(\theta(f_1f_2))=\Phi_{\hspace{-1.1pt}X}(g_1g_2),\vspace{1pt}$ and consequently, $\theta(f_1 f_2)=g_1 g_2=\theta(f_1)\hspace{.4pt}\theta(f_2).$ Since $\theta$ is surjective, $\theta(\id\hspace{.4pt}_I)=\id_X.$
The proof of the theorem is completed.

\medskip

\noindent{\sc Remark.} If $\cC$ is a left triangulated subcategory of a triangulated category, then $\cC\subseteq \cC[1]$. In particular, Theorem \ref{existence-3} covers the essential part of Krause's result stated in \cite[(2.2)]{Kra2}, where the isomorphism $\End(X)\cong \End_{\End(Z)}(I)$ is only verified to be an abelian group isomorphism.

%

\medskip

In a dual fashion, we may establish the following statement.

\begin{Theo}\label{existence-2}

Let $\mathcal{C}$ be an extension-closed subcategory of a triangulated category, and let $X\in \cC$ be strongly indecomposable with a non-zero monomorphism \vspace{-8pt}
$$\Psi :  \Hom_\mathcal{\hspace{.4pt}C}(Z, -) \hspace{-2pt}\to \hspace{-2pt} \Hom_{\hspace{1pt}\overline{\End\hspace{-.5pt}}\hspace{.5pt}(\hspace{-1pt}X\hspace{-1pt})} \hspace{-1pt} (\hspace{.4pt}\oHom_\mathcal{\hspace{.6pt}C }(-, X), I) \vspace{-2pt}$$ where $Z\in \cC\cap \cC[-1]$ and $I$ is an injective envelope of the left $\overline{\End\hspace{-.5pt}}(X)$-module $S_X$. \vspace{-2pt} If $\hspace{.8pt}\Psi_{\hspace{-1.1pt}Z}$ is bijective, then $\cC$ has an almost split sequence $\hspace{-2pt}\xymatrixcolsep{18pt}\xymatrix{X \ar[r] & Y \ar[r] & Z[1].}$

\end{Theo}

\section{Auslander-Reiten functors}

\medskip

\noindent The objective of this section is to study the existence of almost split sequences in a Hom-reflexive tri-exact $R$-category, where $R$ is a commutative ring. Our main results will relate the global existence of almost split sequences in the category to the existence of an Auslander-Reiten functor, which is a generalization of an Auslander-Reiten duality considered in \cite{LeZ}; and in the left or right triangulated case, to the existence of a Serre functor with a proper image, which differs slightly from the classical notion of a Serre functor defined in \cite[(I.1)]{RvdB}; see also \cite{BoK}.

\medskip

Throughout this section, $\cC$ will stand for a tri-exact $R$-category, say an extension-closed subcategory of a triangulated $R$-category $\cA$. We shall say that $\cC$ has {\it almost split sequences on the right} (respectively, {\it left}) if every strongly indecomposable not Ext-projective (respectively, not Ext-injective) object is the ending (respectively, starting) term of an almost split sequence; and that $\cC$ {\it has almost split sequences} if it has almost split sequences on the right and on the left. Recall that the exact functor $D=\Hom_R(-, I\hspace{-1pt}_R): \Mod R\to \Mod R$ restricts to a duality $D: {\rm RMod}\hspace{.5pt}R \to {\rm RMod}\hspace{.5pt}R$, where $I\hspace{-1pt}_R$ is the minimal injective co-generator for $\Mod R$, whereas ${\rm RMod}\hspace{.5pt}R$ is the category of reflexive $R$-modules. We shall say that $\cC$ is {\it Hom-reflexive} (respectively, {\it Hom-finite}) if $\Hom\hspace{.5pt}_\mathcal{C}\hspace{-1pt}(X, Y)$ is reflexive (respectively, of finite length) over $R$, for all $X, Y\in \cC$; and {\it Ext-reflexive} if $\Ext^1_\mathcal{C}\hspace{-1pt}(X, Y)$ is reflexive over $R$, for all $X, Y\in \cC$. Hom-finite $R$-categories are Hom-reflexive; see (\ref{alg-ref-mod}), and the converse is not true. The following statement is important for our purpose; compare \cite[(2.4)]{LNP}.


\begin{Lemma}\label{lin-alg}

Let $R$ be a commutative ring, and let $M, N\in \Mod R$ defining a non-degenerate $R$-bilinear form
$<\hspace{-3pt}-,-\hspace{-3pt}>: M\times N \to I_{\hspace{-1pt}R}.$ If $M$ or $N$ is reflexive, then both $M$ and $N$ are reflexive with $R$-linear isomorphisms $\phi_M: M\to DN: u\mapsto <\hspace{-2.5pt}u, -\hspace{-4pt}>$ and $\psi\hspace{-.5pt}_N: N\to DM: v\mapsto <\hspace{-2.5pt}-, v\hspace{-2.5pt}>$.

\end{Lemma}

\noindent{\it Proof.} By the hypothesis, $\phi_M$ and $\psi\hspace{-.5pt}_N$ are monomorphisms. Consider the canonical monomorphisms $\sigma_{\hspace{-1pt}_M}: M\to D^2M$ and $\sigma_{\hspace{-1pt}_N}: N\to D^2N$. It is easy to verify that $\phi_M=D(\psi\hspace{-.5pt}_N)\circ \sigma_{\hspace{-1pt}_M}$ and $\psi\hspace{-.5pt}_N=D(\phi_M)\circ \sigma_{\hspace{-1pt}_N}$, where $D(\phi_M)$ and $D(\psi\hspace{-.5pt}_N)$ are surjective. If $M$ is reflexive, so are $DM$ and $N$; see (\ref{alg-ref-mod}). In particular, $\sigma_{\hspace{-1pt}_M}$ and $\sigma_{\hspace{-1pt}_N}$ are surjective, so are $\phi_M$ and $\psi\hspace{-.5pt}_N$.
The proof of the lemma is completed.

\medskip

We shall first strengthen the results on the existence of an individual almost split sequence under the Hom-reflexive setting. The following preparatory result is well-known under some classical settings; see, for example, \cite{GR, LeZ, RvdB}.


\begin{Lemma}\label{bilin-forms}

Let $\mathcal{C}$ be a Hom-reflexive tri-exact $R$-category. Consider an almost-zero extension $\delta\in \Ext^1_\sC(Z, X)$ and a linear form $\theta\in D\Ext^1_\sC(Z, X)$ such that $\theta(\delta)\ne 0$. Given any object $L\in \cC$, there exist natural $R$-linear isomorphisms \vspace{-.5pt}
$$\Oa_{L, X}: \oHom_{\sC}(L, X) \to D\Ext^1_{\sC}(Z, L):
\bar{\hspace{-1.5pt}g\hspace{-.5pt}} \mapsto \theta \circ \Ext_{\sC}^1(Z, \, \bar{\hspace{-1.5pt}g\hspace{-.5pt}})\vspace{-8pt}$$
and \vspace{-1.5pt}
$$\Ta_{Z,L}: \uHom_{\sC}(Z, L) \to D\Ext^1_{\sC}(L, X): \underline{f\hspace{-2pt}}\, \mapsto \theta \circ
\Ext_\cC^1(\hspace{.4pt}\underline{f\hspace{-2pt}}\hspace{2pt}, X).
\vspace{3pt}$$

\end{Lemma}

\noindent{\it Proof.} Given $L\in \cC$, by Lemma \ref{art-1}, we obtain two non-degenerate $R$-bilinear forms \vspace{-6pt}
$$<\hspace{-3pt}-, -\hspace{-3.2pt}>\hspace{-2pt}_L:  \oHom_\mathcal{\hspace{.5pt}C}(L, X) \times \Ext_{\mathcal{C}}^1(Z, L)\to I_R:
(\hspace{1.5pt}\bar{\hspace{-1.5pt}g}, \zeta) \mapsto \theta (\hspace{1.5pt}\bar{\hspace{-1.5pt}g}\cdot \zeta)\vspace{-5pt}$$
and \vspace{-2pt}
$${_L}\hspace{-3pt}<\hspace{-3pt}-, -\hspace{-3.2pt}>: \Ext^1_{\sC}(L, X)  \times \uHom_{\sC}(Z, L)  \to I_R: (\zeta, \,\underline{f\hspace{-2pt}}\hspace{2pt})
\to \theta (\zeta \cdot \underline{f\hspace{-2pt}}\hspace{2pt}).$$

Since $\cC$ is Hom-reflexive, by Lemma \ref{alg-ref-mod}, $\oHom_\mathcal{\hspace{.5pt}C}(L, X)$ and $\uHom_{\sC}(Z, L)$ are reflexive $R$-modules. By Lemma \ref{lin-alg}, we obtain two isomorphisms $\Oa_{L, X}$ and $\Ta_{Z,L}$ as stated in the lemma, which are clearly natural in $L$; see (\ref{Ext-bimod}). The proof of the lemma is completed.

\medskip

The following result improves Theorem \ref{existence} under the Hom-reflexive setting.


\begin{Theo}\label{art-6}

Let $\mathcal{C}$ be a Hom-reflexive tri-exact $R$-category with $X, Z\in \cC$ strongly indecomposable. The following statements are equivalent.

\vspace{-2pt}

\begin{enumerate}[$(1)$]

\item There exists an almost split sequence $\hspace{-3pt}\xymatrixcolsep{20pt}\xymatrix{X\ar[r]&Y\ar[r] &Z}\hspace{-3.5pt}$ in $\cC.$

\item There exists a non-zero isomorphism $\Oa_X: \oHom_\mathcal{\hspace{.7pt}C}(-, X)\to D\Ext^1_\mathcal{C}(Z, -)$.

\vspace{.5pt}

\item There exists a non-zero isomorphism $\Ta_Z: \uHom_\mathcal{\hspace{.6pt}C}(Z,-)\to D\Ext^1_\mathcal{C}(-,X).$

\end{enumerate}

\end{Theo}

\noindent{\it Proof.} Given an almost-zero extension $\delta\in\Ext^1_\sC(Z, X)$, we choose $\theta\in D\Ext^1_\sC(Z, X)$ such that $\theta(\delta)\ne 0$. By Lemma \ref{bilin-forms}, we see that $\uHom_\mathcal{\hspace{.6pt}C}(Z,-)\cong D\Ext^1_\mathcal{C}(-,X)$ and $\oHom_\mathcal{\hspace{.7pt}C}(-, X)\cong D\Ext^1_\mathcal{C}(Z, -)$. Thus, Statement (1) implies Statements (2) and (3).

\vspace{1.5pt}

Let now $\Oa_X: \oHom_\mathcal{\hspace{.7pt}C}(-, X)\to D\Ext^1_\mathcal{C}(Z, -)$ be a nonzero isomorphism. In particular, $Z$ is not Ext-projective, and hence, we have a nonzero canonical algebra homomorphism $R\to \underline{\End}(Z)$. Moreover, since $\mathcal{C}$ is Hom-reflexive,  we obtain an isomorphism $\Psi_X: \Ext^1_\mathcal{C}(Z, -)\to D\oHom_\mathcal{\hspace{.7pt}C}(-, X).$
By Theorem \ref{existence}(2), we obtain an almost split sequence as stated in Statement (1). Similarly, we may show that Statement (3) implies Statement (1). The proof of the theorem is completed.

\medskip

\noindent{\sc Remark.} In case $R$ is artinian, Theorem \ref{art-6} is known for an Ext-finite abelian $R$-category; see \cite{GR,LeZ} and for a Hom-finite exact $R$-category; see \cite{LNP}. 

\medskip

We shall weaken the condition that both $X$ and $Z$ are strongly indecomposable stated in Theorem \ref{art-6} in some special cases as below.

\begin{Theo} \label{existence-7}

Let $\mathcal{C}$ be a Hom-reflexive extension-closed subcategory of a triangulated $R$-category.

\vspace{-1.5pt}

\begin{enumerate}[$(1)$]

\item If $X\in \cC$ is strongly indecomposable and $Z\in \cC\cap \cC[1]$, \vspace{-.5pt} then $\cC$ has an almost split sequence $\xymatrixcolsep{18pt}\xymatrix{\hspace{-3pt}X \ar[r] & Y \ar[r] & Z}\hspace{-3pt}$ if and only if $\oHom_\mathcal{\hspace{.6pt}C\hspace{-.8pt}}(-, X)\cong D\Ext^1_\mathcal{\hspace{.6pt}C\hspace{-.5pt}}(Z, -)\ne 0$.

\vspace{1.5pt}

\item If $Z\in \cC$ is strongly indecomposable and $X\in \cC\cap \cC[-1]$, \vspace{-1pt} then $\cC$ has an almost split sequence $\xymatrixcolsep{18pt}\xymatrix{\hspace{-3pt}X \ar[r] & Y \ar[r] & Z}\hspace{-3pt}$ if and only if
$\uHom_\mathcal{\hspace{.8pt}C\hspace{-.8pt}}(Z, -)\cong D \Ext^1_\mathcal{\hspace{.6pt}C\hspace{-.5pt}}(-,X)\ne 0.$

\end{enumerate}

\end{Theo}

\noindent{\it Proof.} We shall only prove the sufficiency of Statement (1). Let $X\in \cC$ be strongly indecomposable and $\oHom_\mathcal{\hspace{.6pt}C\hspace{-.8pt}}(-, X)\cong D\Ext^1_\mathcal{\hspace{.6pt}C\hspace{-.5pt}}(Z, -)\ne 0$, where $Z=M[1]$ for some $M\in \cC$. It suffices to show that $\End(M)$ is local. Since $\Ext^1_\mathcal{C}(Z, -)\cong \Hom_\mathcal{\hspace{.4pt}C}(M, -),$ we obtain an isomorphism $\Psi_X: \oHom_\mathcal{\hspace{.5pt}C}(-, X) \to D\Hom_\mathcal{\hspace{.4pt}C}(M, -)$.\hspace{-4pt} In particular, $\Psi_{X,X}: \overline{\End}(X)\to D\Hom_{\hspace{.6pt}\mathcal{C}}(M, X)$ is a right $\End(X)$-linear isomorphism and $\Psi_{M,X}: \oHom_\sC(M,X)\to D\End(M)$ is a right $\End(M)$-linear isomorphism. Since $M\in \cC[-1]$, by Lemma \ref{left-tri-ideal}(2), we have $\oHom_\mathcal{\hspace{.6pt}C}(M, X)=\Hom_\mathcal{\hspace{.5pt}C}(M, X)$, which is a left $\overline{\End}(X)$-module. Since $\overline{\End\hspace{-1.5pt}}\hspace{1.5pt}(X)$ is reflexive, we obtain a commutative diagram
$$\xymatrixcolsep{38pt}\xymatrixrowsep{16pt}\xymatrix{
\End(M) \ar[d]_\theta \ar[r]^-{\sigma} & D^2\End(M) \ar[d]^-{D(\Psi\hspace{-1pt}_{M,X})}\\
\overline{\End\hspace{-.5pt}}\hspace{.5pt}(X) \ar[r]^-{\Psi_{X,X}} & D\Hom_{\hspace{.6pt}\mathcal{C}}(M, X)
}$$ of $R$-linear isomorphisms, where $\sigma$ is the canonical isomorphism. It remains to show that $\theta$ is an algebra homomorphism. Indeed, we fix arbitrarily a morphism $u\in \Hom_{\hspace{.4pt}\mathcal{C}}(M, X)$. Given any $f\in \End(M)$, in view of the above commutative diagram, we obtain an equation
$$
(1) \hspace{50pt} \Psi\hspace{-1pt}_{X,X}(\theta(f))(u)=\Psi\hspace{-1pt}_{M,X}(u)(f). \hspace{125pt}
$$
On the other hand, consider the commutative diagram
$$\xymatrixrowsep{18pt}\xymatrix{\overline{\End}(X) \ar[r]^-{\Psi_{\hspace{-.5pt}X,X}} \ar[d]_{\Hom_\cC(\bar{u},X)} & D\Hom_\mathcal{\hspace{.8pt}C}(M, X)
\ar[d]^{D\Hom_{\hspace{.5pt}\cC}(M, u)}\\
\oHom_\mathcal{\hspace{1pt}C}(M, X) \ar[r]^{\Psi_{\hspace{-.5pt}M,X}} & D\End(M).}$$
Given any $v\in \End(X)$, we obtain an equation
$$(2) \hspace{50pt} \Psi_{\hspace{-1.1pt}M, X}(\bar{v} \bar{u})(f) =  \Psi_{X,\hspace{-1.5pt}X}(\bar{v})(u f). \hspace{130pt}$$
Let now $f, g\in \End(M)$. Since $\Psi_{X,X}$ is also right $\overline{\End}(X)$-linear, we deduce from the equations (1) and (2) that
$$
\Psi_{X, X}(\theta(f)\theta (g))(u)= \Psi_{X, X}(\theta(f))(\theta (g)u) = \Psi\hspace{-1pt}_{M,X}(\theta (g) \bar{u} )(f) = \Psi_{X,\hspace{-1.5pt}X}(\theta (g))(u f).
$$
Since $\Psi\hspace{-1pt}_{M,X}$ is right $\End(M)$-linear, we deduce from the equation (1) that
$$
\Psi_{X,\hspace{-1.5pt}X}(\theta (g))(u f) =\Psi\hspace{-1pt}_{M,X}(u f)(g)=\Psi\hspace{-1pt}_{M,X}(u)(fg)
=\Psi\hspace{-1pt}_{X,X}(\theta(fg))(u).$$
This yields \vspace{1pt} $\Psi_{X, X}(\theta(f)\theta (g)) =\Psi\hspace{-1pt}_{X,X}(\theta(fg))$, and hence, $\theta(f g )=\theta(f) \hspace{.4pt} \theta(g)$. 
The proof of the theorem is completed.

%


\begin{Cor}\label{ARS-LRT}

Let $\cC$ be a Hom-reflexive extension-closed subcategory of a triangulated $R$-category. \vspace{-1pt}

\begin{enumerate}[$(1)$]

\item If $\cC$ is left triangulated with $X\in \cC$ strongly indecomposable, then $X$ is the starting term of an almost split sequence if and only if the functor $D\oHom_{\hspace{.8pt}\cC}(-, X)$ is representable by a nonzero object in $\cC[-1]$.

\item If $\cC$ is right triangulated with $Z\in \cC$ strongly indecomposable, then $Z$ is the ending term of an almost split sequence if and only if the functor $D\uHom_{\hspace{.6pt}\cC}(Z, -)$ is representable by a nonzero object in $\cC[1]$.

\end{enumerate}

\end{Cor}

\medskip

In order to study the global existence of almost split sequences in $\cC$, we need to generalize the classical notion of an Auslander-Reiten duality; see \cite{ARS, LeZ}.


\begin{Defn}\label{ARF}

Let $\cC$ be a tri-exact $R$-category.

\begin{enumerate}[$(1)$]

\item A {\it right Auslander-Reiten functor} for $\cC$ is a functor $\tau: \uC \to \oC \hspace{2pt}$ with binatural $R$-linear isomorphisms $\Ta_{X, Y} : \uHom_{\bC}(X, Y) \to D\Ext^1_{\sC}(X, \tau Y),$ where $X, Y\in \cC.$

\vspace{1pt}

\item A {\it left Auslander-Reiten functor} for $\cC$ is a functor $\tau^-: \oC\to \uC\hspace{.5pt}$ with binatural $R$-linear isomorphisms $\Oa_{X, Y}: \oHom_{\bC}(X, Y) \to D\Ext^1_{\sC}(\tau^-Y, X),$ where $X, Y\in \cC.$

%

\end{enumerate}

\end{Defn}

\medskip

The following statement collects some properties of an Auslander-Reiten functor.


\begin{Prop}\label{ARF-prop}

Let $\cC$ be a Hom-reflexive tri-exact $R$-category. Then, a right $($or left$)$ Auslander-Reiten functor for $\cC$ is faithful. If $\cC$ is, in addition, right $($or left$)$ triangulated, then a right $($or left$)$ Auslander-Reiten functor is fully faithful.

\end{Prop}

\noindent{\it Proof.} We shall only consider a right Auslander-Reiten functor $\tau: \uC\to  \oC$\vspace{1pt}  with binatural $R$-linear isomorphisms $\Ta_{X, Y}: \uHom_\bC(X, Y)\to D\Ext^1_{\sC}(Y, \tau X)$. Given a morphism $u: X\to Y$, considering the commutative diagram \vspace{-2pt}
$$\xymatrixcolsep{30pt}\xymatrixrowsep{16pt}\xymatrix{
\underline{\End\hspace{-1.5pt}}\hspace{1.5pt}(Y) \ar[d]_{\uHom_{\hspace{.5pt}\mathcal{C}}(\underline{u\hspace{-1.5pt}}\hspace{2pt}, Y)} \ar[r]^-{\Ta_{Y, Y}} & D\Ext^1_\sC(Y, \tau Y)
\ar[d]^{D\Ext^1_\sC(Y, \tau (\underline{u\hspace{-1.5pt}}\hspace{2pt}))} \\
\uHom_\bC(X, Y) \ar[r]^-{\Ta_{X, Y}} & D\Ext^1_\sC(Y, \tau X),}$$ we obtain \vspace{1pt}
$\Ta_{X, Y}(\hspace{.5pt}\underline{u\hspace{-1.5pt}}\hspace{2pt})=\Ta_{Y, Y}(\hspace{.5pt}\underline{\id}\hspace{.8pt}_Y) \circ \Ext^1_\sC(Y, \tau(\underline{u\hspace{-1.5pt}}\hspace{2pt})).$ If $\tau(\hspace{.5pt}\underline{u\hspace{-1.5pt}}\hspace{2pt})=0$, then $\Ta_{X, Y}(\hspace{.5pt}\underline{u\hspace{-1.5pt}}\hspace{2pt})=0$, and hence, $\underline{u\hspace{-1.5pt}}\hspace{2pt}=0$. That is, $\tau$ is faithful. Next, assume that $\cC$ is a Hom-reflexive right triangulated subcategory of a triangulated $R$-category. By Lemmas \ref{left-tri} and \ref{left-tri-ideal}(2), $\cC[1]\subseteq \cC=\oC$. This yields a functor $F=[1]\circ \tau: \uC\to \cC$ with isomorphisms $\Phi_{X,Y}: \uHom_\bC(X, Y)\to D\Hom_\sC(Y, FX),$ which are binatural in $X$ and $Y$. Given $f\in \Hom_\sC(X, Y)$ and $g\in \Hom_\sC(Y, FX)$, considering the commutative diagram \vspace{-2pt}
$$\xymatrixcolsep{50pt}\xymatrixrowsep{16pt}\xymatrix{
\uHom_\bC(Y, FX) \ar[d]_{\Phi_{Y, FX}} & \uHom_\bC(Y, Y) \ar[l]_-{\uHom(Y, \hspace{1pt} \underline{g\hspace{-2pt}}\hspace{2pt})} \ar[d]^-{\Phi_{Y,Y}} \ar[r]^{\uHom(\hspace{.5pt}\underline{f\hspace{-2pt}}\hspace{2pt}, Y)}    &   \uHom_\bC(X, Y) \ar[d]^{\Phi_{X,Y}}  \\
D\Hom_{\sC}(FX, FY) & D\Hom_{\sC}(Y, FY) \ar[l]_-{D\Hom(g, FY)} \ar[r]^-{D\Hom(Y, F(\underline{f\hspace{-2pt}}\hspace{2pt}))}  & D\Hom_{\sC}(Y, FX),}$$
we obtain the following equations
$$(*)\hspace{50pt} \Phi_{Y,FX}(\hspace{.4pt}\underline{g\hspace{-2pt}}\hspace{2pt})(F(\hspace{.5pt}\underline{f\hspace{-2pt}}\hspace{2pt}))=\Phi_{Y,Y}(\underline{\id}\hspace{1pt}_Y)(F(\underline{f\hspace{-2pt}}\hspace{2pt}) g)=\Phi_{X,Y}(\underline{f\hspace{-2pt}}\hspace{2pt})(g).\hspace{40pt}$$

Since $FX\in \cC[1]$, by Lemma \ref{left-tri-ideal}(1), $\Hom_{\sC}(Y, FX)=\uHom\hspace{.6pt}_{\cC}(Y, FX)$. Thus, we obtain an isomorphism $\Phi_{Y, FX}: \Hom_{\sC}(Y, FX)\to D\Hom_\sC(FX,FY),$ which makes \vspace{-2pt}
$$\xymatrixrowsep{16pt}\xymatrix{
\uHom_\sC(X, Y) \ar[r]^-{\Phi_{X,Y}} \ar[d]_-{F} & D\Hom_{\sC}(Y, FX)  \\
\Hom_\sC(FX, FY) \ar[r]^-{\sigma} & D^2\Hom_{\sC}(FX, FY) \ar[u]_-{D(\Phi_{Y, FX})}
}$$ commute, where $\sigma$ is the canonical isomorphism. Indeed, for any $f\in \Hom_\bC(X, Y)$ and $g\in \Hom_\sC(Y, FX)$, in view of the equations in $(*)$, we see that

$$D(\Phi_{Y,FX}) (\sigma(F(\underline{f\hspace{-2pt}}\hspace{2pt})))(g) = \sigma(F(\underline{f\hspace{-2pt}}\hspace{2pt})) (\Phi_{Y,FX}(\hspace{.4pt}\underline{g\hspace{-2pt}}\hspace{2pt})) =
\Phi_{Y,FX}(\hspace{.5pt}\underline{g\hspace{-2pt}}\hspace{2pt})(F(\hspace{.4pt}\underline{f\hspace{-2pt}}\hspace{2pt}))
=\Phi_{X,Y}(\hspace{.4pt}\underline{f\hspace{-2pt}}\hspace{2pt})(g). \vspace{-2pt}$$
As a consequence, $F: \uHom_\sC(X, Y)\to \Hom_\sC(FX, FY)$ is an isomorphism. That is, $F$ is fully faithful, and so is $\tau $. The proof of the proposition is completed.

\medskip

We are now able to relate the existence of almost split sequences to the existence of an Auslander-Reiten functor.


\begin{Theo}\label{ARS-ARF}

Let $\cC$ be a Hom-reflexive Krull-Schmidt tri-exact $R$-category.

\vspace{-2pt}

\begin{enumerate}[$(1)$]

\item There exist almost split sequences on the right $($respectively, left$\hspace{.5pt})$ in $\cC$ if and only if it admits a full right $($respectively, left$\hspace{.5pt})$ Auslander\--Reiten functor.


\vspace{1pt}

\item There exist almost split sequences in $\cC$ if and only if it admits a right Auslander-Reiten equivalence,
or equivalently, a left Auslander-Reiten equiva\-lence$\,;$ and in this case, $\cC$ is Ext-reflexive.

\end{enumerate}

\end{Theo}

\noindent{\it Proof.} We shall prove the theorem only for right Auslander-Reiten functors. Consider a full right Auslander-Reiten functor $\tau: \uC\to  \oC$. Let $Z\in \cC$ be indecomposable but not Ext-projective. By defi\-nition, there exist $R$-linear isomorphisms $\Ta_{Z,Y}: \uHom_{\sC}(Z,Y)\to D\Ext^1_{\sC}(Y, \tau Z)$, which is natural in $Y$. Therefore, $\uHom_{\bC}(Z, - )\cong D\Ext^1_{\sC}(-, \tau Z).$ By Proposition \ref{ARF-prop}, $\overline{\End\hspace{-1.5pt}}\hspace{1.5pt}(\tau Z)\cong \underline{\End\hspace{-1pt}}\hspace{1.5pt}(Z)$, which is local. By Lemma \ref{Lift-stab-obj}(1), $\Ext_\cC^1(-, \tau Z) \cong \Ext_\cC^1(-, X)$ for some indecomposable $X\in \cC$. Then, $\uHom_{\bC}(Z, -)\cong D\Ext^1_{\sC}(-, X).$
By Theorem \ref{art-6}, there exists an almost split sequence $\hspace{-3pt}\xymatrixcolsep{20pt}\xymatrix{X\ar[r]& L \ar[r] & Z}$ in $\cC$.

Assume that $\cC$ has almost split sequences on the right. For each indecomposable and not Ext-projective object $M\in \cC$, we fix an indecomposable object $\tau M\in \cC$, an almost-zero extension $\delta_M \in \Ext^1_\sC(M, \tau M)$, and a linear form
$\theta_M\in D\Ext^1_\sC(M, \tau M)$ such that $\theta_M(\delta_M)\ne 0$.
Let $X, Y, Z\in \cC$ be indecomposable and not Ext-projective. Considering $\theta_X$ and $\theta_Y$, by Lemma \ref{bilin-forms}, we obtain two $R$-linear isomorphisms
$$\Ta_{X,Y}: \uHom_{\sC}(X, Y) \to D\Ext^1_{\sC}(Y, \tau X): \underline{f\hspace{-2pt}}\, \mapsto \theta_X \circ
\Ext_\cC^1(\,\underline{f\hspace{-2pt}}\hspace{2pt}, X)\vspace{-3pt} \hspace{15pt}$$
and \vspace{-2pt}
$$\Oa_{\tau X, \tau Y}: \oHom_{\sC}(\tau X, \tau Y) \to D\Ext^1_{\sC}(Y, \tau X): \hspace{1.5pt}\bar{\hspace{-1.5pt}g} \mapsto \theta_Y \circ \Ext_{\sC}^1(Y, \, \bar{g}).$$
This yields an $R$-linear isomorphism
$$\tau_{\hspace{-.5pt}_{X, Y}}=\Oa_{\tau X, \tau Y}^{-1} \Ta_{X,Y} \hspace{-2pt}: \uHom_\bC(X, Y)  \to \oHom_\bC(\tau X, \tau Y).$$

\vspace{1pt}

Let $f\in \Hom_{\sC}(X, Y)$ and $\zeta\in \Ext^1_\cC(Y, \tau X)$. Since $\Ta_{X, Y}(\underline{f\hspace{-2pt}}\hspace{2pt})=\Oa_{\tau X, \tau Y}(\tau_{\hspace{-.5pt}_{X, Y}}\hspace{-1.5pt}(\hspace{.5pt}\underline{f\hspace{-2pt}}\hspace{2pt}))$ by definition, we obtain
an equation
$$(*) \hspace{20pt} \theta_X(\zeta \cdot \underline{f\hspace{-2pt}}\hspace{2pt})=\theta_Y(\tau_{\hspace{-.5pt}_{X, Y}}\hspace{-1.5pt}(\hspace{1pt}\underline{f\hspace{-2pt}}\hspace{2pt}) \cdot \zeta ).$$

Let $g\in \Hom_{\sC}(Y, Z)$ and $\delta\in \Ext_\cC^1(Z, \tau X)$. By the above equation, we obtain
$$\theta_X(\delta\cdot (\underline{g\hspace{-2pt}}\,\underline{f\hspace{-2pt}}\hspace{2.5pt}))
=\theta_X((\delta\cdot \underline{g\hspace{-2pt}}\hspace{2.5pt})\cdot \underline{f\hspace{-2pt}}\hspace{2.5pt})
=\theta_Y(\tau_{\hspace{-.5pt}_{X, Y}}\hspace{-1.5pt}(\hspace{1pt}\underline{f\hspace{-2pt}}\hspace{2pt}) \cdot \delta\cdot  \underline{g\hspace{-2pt}}\hspace{2pt})
=\theta_Z((\tau_{\hspace{-.5pt}_{Y, Z}}\hspace{-1.5pt}(\hspace{1pt}\underline{g\hspace{-2pt}}\hspace{2pt}) \, \tau_{\hspace{-.5pt}_{X, Y}}\hspace{-1.5pt}(\hspace{1pt}\underline{f\hspace{-2pt}}\hspace{2pt})) \cdot \delta)$$
That is, $\Ta_{X,Z}(\hspace{.5pt}\underline{g\hspace{-2pt}}\,\underline{f\hspace{-2pt}}\hspace{2.5pt}) (\delta)=\Oa_{X,Z}(\tau_{\hspace{-.5pt}_{Y, Z}}\hspace{-1.5pt}(\hspace{1pt}\underline{g\hspace{-2pt}}\hspace{2pt}) \, \tau_{\hspace{-.5pt}_{X, Y}}\hspace{-1.5pt}(\hspace{1pt}\underline{f\hspace{-2pt}}\hspace{2pt}))(\delta),$ from which we conclude that
$\tau_{\hspace{-.5pt}_{X,Z}}\hspace{-1.5pt}(\hspace{.5pt}\underline{g\hspace{-2pt}}\,\underline{f\hspace{-2pt}}\hspace{2.5pt})
=\tau_{\hspace{-.5pt}_{Y, Z}}\hspace{-1.5pt}(\hspace{1pt}\underline{g\hspace{-2pt}}\hspace{2.5pt}) \, \tau_{\hspace{-.5pt}_{X, Y}}\hspace{-1.5pt}(\hspace{1pt}\underline{f\hspace{-2pt}}\hspace{2pt}).$ \vspace{.5pt} In particular, since $\tau_{\hspace{-.5pt}_{X, X}}$ is bijective, $\tau_{\hspace{-.5pt}_{X, X}}\hspace{-1.5pt}(\hspace{.8pt}\underline{\id}\hspace{1pt}_X)=\bar{\hspace{-.4pt}\id}_{\hspace{1.2pt} \tau \hspace{-1pt} X}$. Since $\uC$ is Krull-Schmidt, by Lemma \ref{Lift-stab-obj}(1), we may extend $\tau$ to a fully faithful functor $\tau: \uC\to \oC\hspace{.5pt}.$ It remains to show that the isomorphism $\Ta_{X,Y}$ is binatural. It is natural in $Y$ by Lemma \ref{bilin-forms}. We claim, for $h\in \Hom_{\sC}(X, Z)$, that the diagram

$$\xymatrixcolsep{28pt}\xymatrixrowsep{16pt}
\xymatrix{\uHom_{\bC}(Z, Y) \ar[r]^-{\it\Theta_{Z, Y}} \ar[d]_{\uHom_{\hspace{.4pt}\cC}(\underline{h\hspace{-1pt}}\hspace{2pt}, Y)} & D\Ext^1_{\sC}(Y, \tau Z) \ar[d]^{D\Ext^1_{\sC}(Y, \hspace{1pt}
\tau_{\hspace{-1pt}_{X, Z}}\hspace{-.8pt}(\hspace{-.5pt}\underline{h\hspace{-1pt}}\hspace{1pt}))} \\
\uHom_{\bC}(X, Y) \ar[r]^-{\it\Theta_{X, Y}} & D\Ext^1_{\sC}(Y, \tau X)
} $$
is commutative. Indeed, for any $u\in \Hom_{\hspace{.5pt}\cC}(Z, Y)$ and $\delta\in \Ext^1_{\sC}(Y, \tau Z)$, we obtain
$$D\Ext^1_{\sC}(Y, \hspace{1pt} \tau_{\hspace{-1pt}_{X, Z}}\hspace{-1pt}(\hspace{-.5pt}\underline{h\hspace{-1pt}}\hspace{1pt}))
(\Ta_{Z, Y}(\underline{u\hspace{-1pt}}\hspace{1pt}))(\delta)=\Ta_{Z, Y}(\underline{u\hspace{-1pt}}\hspace{1pt})
(\tau_{\hspace{-1pt}_{X, Z}}\hspace{-1pt}(\hspace{-.5pt}\underline{h\hspace{-1pt}}\hspace{1pt})\cdot \delta)=
\theta_Z(\tau_{\hspace{-1pt}_{X, Z}}\hspace{-1pt}(\hspace{-.5pt}\underline{h\hspace{-1pt}}\hspace{1pt})\cdot \delta\cdot \underline{u\hspace{-1pt}}\hspace{1pt}).\hspace{12pt}$$
On the other hand, applying the definition of $\Ta_{X\hspace{-1pt}, Y}$ and the equation $(*)$, we obtain
$$
\Ta_{X\hspace{-1pt}, Y}(\underline{u\hspace{-1pt}} \, \underline{h\hspace{-1pt}}\hspace{1pt})(\delta)=
\theta_X( (\delta \cdot \underline{u\hspace{-1pt}} ) \cdot \underline{h\hspace{-1pt}}\hspace{1pt})
=\theta_Z(\tau_{\hspace{-1pt}_{X, Z}}\hspace{-1pt}(\hspace{-.5pt}\underline{h\hspace{-1pt}}\hspace{1pt})\cdot \delta\cdot \underline{u\hspace{-1pt}}\hspace{1pt}).$$ This shows that $\Ta_{X,Y}$ is natural in $X$. This establishes Statement (1).

Next, assume that the right Auslander-Reiten functor $\tau: \uC\to \oC$ is an equivlence. Since $\tau$ is dense, $\cC$ has almost split sequences on the left.
Let $X, Y\in \cC$ be indecomposable. If $Y$ is Ext-injective, then $\Ext^1_\sC(X, Y)=0$. Otherwise, we obtain $Y\cong \tau Z$, where $Z\in \cC$ is indecomposable and not Ext-projective. In this case,
$D\Ext^1_\sC(X, Y)\cong D\Ext^1_\sC(X, \tau Z)\cong \uHom_\sC(Z, X)$. Since $\Hom_\sC(Z, X)$ is reflexive, by Lemma \ref{alg-ref-mod}, so are $\uHom_\sC(Z, X)$ and $\Ext^1_\sC(X, Y)$. Being Krull-Schmidt, $\cC$ is Ext-reflexive. The proof of the theorem is completed.

\medskip

\noindent{\sc Remark.} Theorem \ref{ARS-ARF}(2) generalizes Lenzing and Zuazua's result stated in \cite[(1.1)]{LeZ} for Ext-finite abelian categories over a commutative artinian ring.

\medskip

\noindent{\sc Example.} Let $Q$ be a quiver of type $\mathbb{A}_\infty$ with a unique source vertex. The category of finitely presented representations of $Q$ over a field $k$ is a Hom-finite abelian $k$-category, which admits a right Auslander-Reiten functor, but no left Auslander-Reiten functor; see \cite[(1.15), (3.7)]{BLP}.

\medskip

Finally, we shall specialize to left or right triangulated $R$-categories. For this purpose, we shall modify the classical notion of a Serre functor; see \cite[(I.1)]{RvdB}.


\begin{Defn}\label{SF}

Let $\cC$ be a tri-exact $R$-category.

\begin{enumerate}[$(1)$]

\item A {\it left Serre functor} for $\cC$ is a functor $\mathbb{S}: \oC\to \uC\hspace{.5pt}$ with binatural $R$-linear isomorphisms $\Phi_{X, Y}: \oHom_{\bC}(X, Y) \to D\uHom_{\bC}(\mathbb{S}Y, X),$ with $X, Y\in \cC.$

\vspace{1pt}

\item A {\it right Serre functor} for $\cC$ is a functor $\mathbb{S}: \uC \to \oC \hspace{2pt}$ with binatural $R$-linear isomorphisms $\Psi_{X, Y} : \uHom_{\bC}(X, Y) \to D\oHom_{\bC}(Y, \mathbb{S}X)$, with $X, Y\in \cC.$

%

\end{enumerate}

\end{Defn}

\medskip

\noindent{\sc Remark.} In case $\cC$ is a triangulated $R$-category, by Lemma \ref{left-tri-ideal},
our left or right Serre functors coincide with those given by Reiten and van den Bergh in \cite[(I.1)]{RvdB}.


\begin{Theo}\label{AR-Serre}

Let $\mathcal{C}$ be a Hom-reflexive Krull-Schmidt tri-exact $R$-category.

\vspace{-3pt}

\begin{enumerate}[$(1)$]

\item If $\cC$ is right triangulated, then it has almost split sequences on the right if and only if it admits a right Auslander-Reiten functor, or equivalently, a right Serre functor whose image lies in $\cC[1].$

\vspace{1pt}

\item If $\cC$ is left triangulated, then it has almost split sequences on the left if and only if it admits a left Auslander-Reiten functor, or equivalently, a left Serre functor whose image lies in $\cC[-1].$

\vspace{1pt}

\item If $\cC$ is triangulated, then it has almost split sequences if and only if it admits a right or left Serre equivalence, or equivalently, a right or left Serre equivalence.

\end{enumerate}

\end{Theo}

\noindent{\it Proof.} Since Statement (3) is an immediate consequence of Statements (1) and (2), we shall only prove Statement (1).
Assume that $\cC$ is a right triangulated subcategory of a triangulated $R$-category. The first equivalence stated in Statement (1) follows from Proposition \ref{ARF-prop} and Theorem \ref{ARS-ARF}(2). For the second equivalence, observe that $\cC[1]\subseteq \cC=\oC$; see (\ref{left-tri}) and (\ref{left-tri-ideal}). Let $\tau: \uC \to \cC$ be a right Auslander-Reiten functor with binatural isomorphisms $\Ta_{X,Y}\hspace{-1.5pt}:\hspace{-.5pt} \uHom_{\sC}(X,Y)\hspace{-1.5pt}\to \hspace{-2pt} D\Ext^1_\sC(Y, \tau X)$. Since $\Ext^1_\sC(Y, \tau X)=\Hom_\sC(Y, (\tau X)[1])$, we see that $\mathbb{S}=[1]\circ \tau: \uC\to \cC$ is a right Serre functor, whose image lies in $\cC[1]$. Conversely, let $\mathbb{S}: \uC \to \cC$ be a right Serre functor with binatural isomorphisms $\Psi_{X,Y}: \uHom_{\sC}(X,Y)\to D\Hom_\sC(Y, \mathbb{S} X)$. Suppose that $\mathbb{S}(\uC)\subseteq \cC[1]$. Then, $\Hom_\sC(Y, \mathbb{S} X)= \Ext^1_\sC(Y, (\mathbb{S} X)[-1])$, for $X, Y\in \cC$. Thus, $\tau=[-1]\circ \mathbb{S}: \uC\to \cC$ is a right Auslander-Reiten functor. The proof of the theorem is completed.

\medskip

\noindent{\sc Remark.} Theorem \ref{AR-Serre} generalizes Reiten and van den Bergh's result stated in \cite[(I.2.4)]{RvdB} for a Hom-finite triangulated category over a field.

\medskip

\noindent{\sc Example.} Let $\La=kQ/J$ be a strongly locally finite dimensional algebra over a field $k$, where $Q$ is a locally finite quiver without infinite paths and $J$ is locally admissible. Since all modules in ${\rm mod}^{\hspace{.3pt}b\hspace{-2.5pt}}\La$ have finite projective and finite injective dimension, we will see from Theorem \ref{ART-refl} that $D^b({\rm mod}^{\hspace{.3pt}b\hspace{-2.5pt}}\La)$ has almost split sequences. Further, for each $n\in \Z$, the right triangulated category $D^{\le n}({\rm mod}^{\hspace{.3pt}b\hspace{-2.5pt}}\La)$ has almost sequences on the left, and the left triangulated category $D^{\ge n}({\rm mod}^{\hspace{.3pt}b\hspace{-2.5pt}}\La)$ has almost sequences on the right. More examples can be found at the end of this paper.

\section{Almost split triangles in derived categories}

\medskip

\noindent The main objective of this section is to study almost split triangles in the derived categories of an abelian category with enough projective objects and enough injective objects. Our results are applicable to the derived categories of modules categories over an algebra with a unity or a locally finite dimension algebra given by a quiver with relations. In particular, they include Happel's result obtained in \cite{Hap1} for the bounded derived category of finite dimensional modules over a finite dimensional modules.
%
%
We shall start with an arbitrary abelian category $\mf A$ and quote the following well-known statement; see, for example, \cite[(10.4.7)]{WEI}.


\begin{Lemma}\label{Hom-iso}

Let $X^\ydt, Y^\ydt$ be complexes over an abelian category $\mf{A}$. If $X^\ydt$ is bounded-above of projective objects or $Y^\ydt$ is bounded-below of injective objects, then there exists an isomorphism
$\mathbb{L}_{X^\ydt,Y^\ydt}:\Hom_{K(\mf{A})}(X^\ydt, Y^\ydt) \to \Hom_{D(\mf{A})}(X^\ydt, Y^\ydt),$ which is induced from the localization functor $\mathbb{L}: K(\mf{A})\to D(\mf{A})$.

\end{Lemma}

\medskip

Let $\mathcal P$ and $\mathcal I$ be strictly additive subcategories of $\mf A$ of projective objects and of injective objects, respectively.
By Lemma \ref{Hom-iso}, we can view $K^b(\mathcal P)$ and $K^b(\mathcal I)$ as full subcategories of $D^b(\mf A)$. A {\it projective resolution} over $\mathcal P$ of a complex $Z^\pdt\in C^-(\mf A)$ is a quasi-isomorphism $s^\sdt: P^\pdt\to Z^\pdt$ with $P^\pdt\in C^-(\mathcal P),$ which is {\it finite} if $P^\pdt\in C^b(\mathcal P).$ Dually, an injective {\it co-resolution} over $\mathcal I$ of a complex $X^\ydt\in C^+(\mf A)$ is a quasi-isomorphism $t^\sdt: X^\ydt \to I^\ydt$ with $I^\ydt\in C^+(\mathcal I),$ which is {\it finite} if $I^\ydt\in C^b(\mathcal I).$


\begin{Theo}\label{ART-nec}

Let $\mf{A}$ be an abelian category such that
$D^*(\mf{A})$ with $*\in \{\emptyset, +, -, b\}$ has an almost split triangle $\hspace{-2pt}\xymatrixcolsep{18pt}\xymatrix{X^\ydt\ar[r] & Y^\ydt \ar[r] & Z^\ydt \ar[r] & X^\ydt[1],}\hspace{-2pt}\vspace{-1pt}$ where $X^\ydt$ is a bounded-below complex and $Z^\ydt$ is a bounded-above complex.

\vspace{-1pt}

\begin{enumerate}[$(1)$]

\item If $Z^\pdt$ admits a projective resolution over a strictly additive subcategory $\mathcal P$ of projective objects of $\mf A$, then it admits a finite projective resolution over  $\mathcal P$.

\vspace{.5pt}

\item If $X^\ydt$ admits an injective co-resolution over a strictly additive subcategory $\mathcal I$ of injective objects of $\mf A$, then it admits a finite injective co-resolution over $\mathcal I$.

\end{enumerate}

\end{Theo}

\noindent{\it Proof.} We shall view $D^*(\mf{A})$ as a full triangulated subcategory of $D(\mf{A})$; see \cite[Chapter III]{Mil}. Let $\mathcal P$ be a strictly additive subcategory of projective objects $\mf A$ with $s^\sdt: P^\pdt \to Z^\pdt$ a quasi-isomorphism with $P^\pdt\in C^-({\mathcal P})$. Write $W^\pdt=X^\cdt[1]$, a complex in $D^+(\mf A)\cap D^*(\mf A)$. Let $n$ be an integer such that $W^i=0$ for all $i < n$. Write $\delta^\pdt: Z^\pdt\to X^\pdt[1]$ for the third morphism in the almost split triangle stated in the theorem. By Lemma \ref{Hom-iso}, $\delta^\pdt \hspace{.5pt} \tilde{s}^\sdt=\tilde{t}^\sdt$ for some complex morphism $t^\pdt: P^\pdt\to W^\pdt$.

Consider the brutal truncation $\kappa_{\ge n}(P^\pdt)$ and the associated canonical morphism $\mu^\ydt: \kappa_{\ge n}(P^\pdt) \to P^\pdt$. Being bounded, $\kappa_{\ge n}(P^\pdt)$ lies in $D^*(\mf A)$. We claim that $\bar \mu^\ydt$ is a retraction in $K(\mf{A})$. Otherwise, by Lemma \ref{Hom-iso}, $\tilde{s}^\sdt \hspace{.5pt} \tilde{\mu}^\sdt$ is not a retraction in $D^*(\mf{A})$, and hence, $\tilde t^\sdt \hspace{.5pt} \tilde \mu^\pdt = \delta^\sdt (\tilde{s}^\sdt \hspace{.5pt} \tilde\mu^\sdt)=0$. By Lemma \ref{Hom-iso}, $\bar t^\sdt \hspace{.5pt} \bar\mu^\sdt=0$. In particular, there exist morphisms $h^i: P^i\to W^{i-1}$ with $i\ge n$ such that $t^i=t^i\mu^i=h^{i+1} d_P^i + d_W^{i+1}h^i,$ \vspace{.5pt}
for all $i\ge n.$ Setting $h^i=0: P^i\to W^{i-1}$ for $i<n$, we obtain $t^i=h^{i+1}d_P^i + d_W^{i+1}h^i,$ for all $i\in \Z$. That is, $\bar t^\sdt=0,$ and hence, $\delta^\pdt=0$, a contradiction. This establishes our claim. In particular, ${\rm H}^i(P^\pdt)=0$ for all $i<n$.

Let $u^\cdt: P^\sdt\to \kappa_{\ge n}(P^\sdt)$ be a complex morphism such that $\bar\mu^\cdt \hspace{.5pt} \bar u^\cdt=\bar{\hspace{-.5pt}\id}_{P^\sdt}$. Then, there exist $f^{i+1}: P^{i+1}\to P^i$ such that $1_{\hspace{-1pt}P^i}-u^i=f^{i+1} d_P^i+d_P^{i-1} f^i$, for $i\in \Z$. In particular, $1_{P^{n-1}}=f^n d_P^{n-1}+d_P^{n-2}f^{n-1}$ and $d_P^{n-1}=d^{n-1}_P f^n d_P^{n-1}.$ Write $d^{n-1}_P= j v,$ where
$v: P^{n-1}\to C$ is the cokernel of $d_P^{n-2}$. Since ${\rm Im}(d_P^{n-2})={\rm Ker}(d_P^{n-1})$, by the Snake Lemma, $j: C\to P^n$ is a monomorphism. Since $ju=juf^nju$, we obtain $1_Q=(uf^n)j$, and hence, $C\in \mathcal{P}$. Thus, the smart truncation $\tau_{\hspace{.5pt}\ge n}(P^\pdt)$ lies in $C^b(\mathcal{P})$. Since ${\rm H}^i(P^\pdt)=0$ for all $i<n$, the canonical projection $p^\sdt: P^\pdt \to \tau_{\hspace{.5pt}\ge n}(P^\sdt)$ is a quasi-isomorphism; see \cite[(III.3.4.2)]{Mil}. Therefore, $Z^\ydt\cong \tau_{\hspace{.5pt}\ge n}(P^\sdt)$ in $D^*(\mf{A})$. By Lemma \ref{Hom-iso}, $\tau_{\hspace{.5pt}\ge n}(P^\pdt)$ is a finite projective resolution of $Z^\ydt$ over $\mathcal P$. Dually, we may establish Statement (2). The proof of the theorem is completed.

\medskip

If $\mf A$ has enough projective (respectively, injective) objects, then every bounded-above (respectively, bounded-below)
complex over $\mf A$ admits a projective resolution (respectively, injective co-resolution); see \cite[(7.5)]{BGKHME}.


\begin{Cor}\label{ART-bd}

Let $\mf A$ be an abelian category \vspace{-2pt} such that $D^b(\mf{A})$ has an almost split triangle $\hspace{-2pt}\xymatrixcolsep{18pt}\xymatrix{X^\pdt\ar[r] & Y^\pdt \ar[r] & Z^\pdt \ar[r] & X^\pdt\hspace{.4pt}[1].}$

\vspace{-6pt}

\begin{enumerate}[$(1)$]

\item If $\mf A$ has enough projective objects, then $Z^\pdt$ has a finite projective resolution.

\vspace{0pt}

\item If $\mf A$ has enough injective objects, then $X^\pdt$ has a finite injective co-resolution.

\end{enumerate}

\end{Cor}

\smallskip

\noindent{\sc Example.} Given any ring $\Sa$, Corollary \ref{ART-bd} applies in $D^b(\Mod \Sa)$; and if $\Sa$ is noetherian, then Corollary \ref{ART-bd}(1) applies in $D^b({\rm mod}^+\hspace{-2pt} \Sa)$.

\medskip

Next, we shall obtain some sufficient conditions for the existence of an almost split triangle in the derived categories of $\mf A$. For this purpose, we need to assume that $\mf A$ is an abelian $R$-category and consider $D={\rm Hom}_R(-, I\hspace{-1.8pt}_R): \Mod R \to \Mod R,$ where $I\hspace{-1.8pt}_R$ is a minimal injective co-generator for $\Mod R$.


\begin{Defn}\label{Naka-Func}

Let $\mf A$ be an abelian $R$-category. Given $\mathcal P$ a strictly additive subcategory of projective objects of $\mf A$,
a functor $\nu: \mathcal{P}\to {\mf A}$ is called a {\it Nakayama functor} if there exist binatural isomorphisms $\beta_{\hspace{-1pt}_{P, X}}:\Hom_{\mathfrak A}(X, \nu P)\to D\Hom_{\mathfrak A}(P, X)$, for all $P\in {\mathcal P}$ and $X\in {\mf A}$.

\end{Defn}

\medskip

\noindent{\sc Remark.} Given a Nakayama functor $\nu: {\mathcal P}\to \mf A$, we see easily that $\nu P$ is an injective object of $\mf A$, for every $P\in \mathcal P$. Hence $\nu \mathcal P$, the image of $\mathcal P$ under $\nu$, is a strictly additive subcategory of injective objects of $\mf A$.

\medskip

As an example, we have the following probably known statement.

\begin{Lemma}\label{Alg-nf} Let $A$ be an $R$-algebra. Then $\nu_{\hspace{-1.5pt}_A}=D\Hom_A(-, A): {\rm proj}\hspace{.5pt}A \to \Mod A$ is a Nakayama functor for $\Mod A$.

\end{Lemma}

\noindent{\it Proof.} Given $P\in {\rm proj}\hspace{.5pt} A$ and $X\in \Mod A$, it is well known; see \cite[(20.10)]{AnF} that there exists a binatural $R$-linear isomoprhism $$\eta_{_{P,X}}: \Hom_A(P, A)\otimes_A X\to \Hom_A(P, X): f\otimes x \mapsto [\hspace{.5pt}u\mapsto f(u) x \hspace{.5pt}].$$
Considering the $R$-$A$-bimodule $\Hom_A(P, A)$ and the adjoint isomorphism, we obtain the following binatural isomorphisms
$$\begin{array}{rcl}
\Hom_R(\Hom_A(P, X), I\hspace{-1.5pt}_R) & \stackrel{\sim}{\longrightarrow} & \Hom_R(\Hom_A(P, A)\otimes_AX, I\hspace{-1.5pt}_R) \vspace{1pt} \\
& \stackrel{\sim}{\longrightarrow} &\Hom_A(X, \Hom_R(\Hom_A(P, A), I\hspace{-1.5pt}_R)).
\end{array}$$
The proof of the lemma is completed.

\medskip

The following statement collects some properties of a Nakayama functor.

\begin{Lemma}\label{NFPro}

Let $\mf A$ be an abelian category with $\mathcal P$ a strictly additive subcategory of projective objects of $\mf A$. Then every Nakayama functor $\nu: {\mathcal P}\to {\mf A}$ is faithful, and it is fully faithful in case $\mathcal P$ is Hom-reflexive over $R$.

\end{Lemma}

\noindent{\it Proof.} Let $\nu: {\mathcal P}\to {\mf A}$ be a Nakayama functor with binatural $R$-linear isomorphisms $\beta_{\hspace{-.8pt}_{P, X}}: \vspace{.5pt} \Hom_{\mathfrak A}(X, \nu P)\to D\Hom_{\mathfrak A}(P, X),$ where $P\in {\mathcal P}$ and $X\in {\mf A}$.
Fix two objects $L, P \in \mathcal P$. Given $f: L\to P$ and $g: P\to \nu L$, \vspace{-3pt} considering the commutative diagram
$$\xymatrixcolsep{26pt}\xymatrixrowsep{16pt}\xymatrix{
\Hom_{\hspace{.5pt}\mathfrak A}(\nu L, \nu P) \ar[d]^{\beta_{P, \nu L}} & \Hom_{\mathfrak A}(\nu L, \nu L) \ar[r]^-{(gf, \nu L)} \ar[d]^-{\beta_{L,\nu L}} \ar[l]_-{(\nu L, \nu f)} &    \Hom_{\mathfrak A}(L, \nu L) \ar[d]^{\beta_{L, L}} & \Hom_{\mathfrak A}(P, \nu L) \ar[d]^{\beta_{L, P}}  \ar[l]_-{(f, \nu L)} \\
D \Hom_{\mathfrak A}(P, \nu L) & D\Hom_{\mathfrak A}(L, \nu L) \ar[r]^-{D(L, gf)} \ar[l]_-{D(f, \nu L)} &
D \Hom_{\mathfrak A}(L, L) & D \Hom_{\mathfrak A}(L, P) \ar[l]_-{D(L, f)},}$$
we obtain the following equations
$$\hspace{30pt} (*) \hspace{30pt} \beta_{P, \nu L}(\nu f)(g)=\beta_{L,\nu L}(\id_{\hspace{.5pt}\nu L})(gf)=\beta_{L, L}(gf)(\id_L)=\beta_{L, P}(g)(f).\hspace{50pt}$$
We claim these equations imply the commutativity of the diagram
$$\xymatrixcolsep{32pt}\xymatrixrowsep{16pt}\xymatrix{
\Hom_{\mathfrak A}(L, P) \ar[r]^-{\sigma} \ar[d]_\nu& D^2 \Hom_{\mathfrak A}(L, P)\ar[d]^{D(\beta_{L,P})}\\
\Hom_{\mathfrak A}(\nu L, \nu P) \ar[r]^-{\beta_{P, \nu L}} & D \Hom_{\mathfrak A}(P, \nu L),}$$ where $\sigma$ is the canonical injection. Indeed, using the equations in $(*)$, we see that
$$D(\beta_{L,P})(\sigma(f))(g)=\sigma(f)(\beta_{L,P}(g))=\beta_{L,P}(g)(f)=\beta_{P, \nu L}(\nu f)(g).$$
As a consequence, $\nu: \Hom_{\mathfrak A}(L, P) \to \Hom_{\mathfrak A}(\nu L, \nu P)$ is a monomorphism, and it is an isomorphism if $\Hom_{\mathfrak A}(L, P)$ is reflexive. The proof of the lemma is completed.

\medskip

\noindent{\sc Remark.} If $\mathcal{P}$ is Hom-reflexive over $R$, then every Nakayama functor $\nu: \mathcal{P}\to \mf A$ co-restricts an equivalence $\nu: \mathcal{P} \to \nu \mathcal{P}$. In this case, we shall always denote by $\nu^{\hspace{.4pt}\mbox{-}}: \nu\mathcal{P} \to \mathcal{P}$ a quasi-inverse of $\nu: \mathcal{P} \to \nu \mathcal{P}$.

\medskip

\noindent{\sc Example.} Let $\La=kQ/I$ be a locally finite dimensional algebra over a field $k$, where $Q$ is locally finite and $J$ is weakly admissible. Then ${\rm proj}\hspace{.5pt}\La$ is Hom-finite over $k$; see \cite[(3.2)]{BHL}, and we have a Nakayama functor $\nu_{\hspace{-1.5pt}_{\mathit\Lambda}}: {\rm proj}\hspace{.5pt}\La \to \Mod \La$,
sending $P_x$ to $I_x$; see \cite[(3.2), (3.6)]{BHL}. This yields an equivalence
$\nu_{\hspace{-1.5pt}_{\mathit\Lambda}}: {\rm proj}\hspace{.5pt}\La \to {\rm inj}\hspace{.5pt}\La$ with a quasi-inverse $\nu^{\mbox{-}}_{\hspace{-1.5pt}_{\mathit\Lambda}}: {\rm inj}\hspace{.5pt}\La \to {\rm proj}\hspace{.5pt}\La$,
sending $I_x$ to $P_x$.

%
%
%
%
%

\medskip

Let $F: \mathcal{A}\to \mathcal{B}$ be a functor between additive categories. Applying $F$ component-wise, one may extend $F$ to a functor  $C(\cA)\to C(\cB)$, sending null-homotopic morphisms to null-homotopic ones and cones to cones. The latter functor induces a triangle-exact functor $K(\cA)\to K(\cB)$; see \cite[(V.1.1.1)]{Mil}. For the simplicity of notation, these functors will be written as $F: C(\cA)\to C(\cB)$ and $F: K(\cA)\to K(\cB)$.


\begin{Prop}\label{Hotopy-fun}

Let $\mf A$ be an abelian $R$-category admitting a Nakayama functor $\nu: {\mathcal P}\to \mf A$, where $\mathcal P$ is a strictly additive subcategory of projective objects of $\mf A$.

\vspace{-2pt}

\begin{enumerate}[$(1)$]

\item The triangle-exact functor $\nu:\hspace{-1pt} K(\mathcal P) \hspace{-1pt}\to\hspace{-1pt} K(\mf A)$ restricts to a triangle-exact functor $\nu :  K^b(\mathcal P) \to  K^b(\nu \mathcal P),$ which is an equivalence if $\mathcal P$ is Hom-reflexive over $R$.

    \vspace{.5pt}

\item Given a complex $X^\ydt$ over $\mf A$ and a bounded complex $P^\pdt$ over $\mathcal P$, we obtain a bina\-tural $R$-linear isomorphism $\tilde{\beta}_{_{\hspace{-1pt}P^\pdt, X^\ydt}}: {\rm Hom}_{D(\mathfrak{A})}(X^\ydt,\nu P^\pdt)\rightarrow
D{\rm Hom}_{D(\mathfrak A)}(P^\pdt, X^\ydt).$

\end{enumerate}

\end{Prop}

\noindent{\it Proof.} If $\mathcal P$ is Hom-reflexive over $R$, then $\nu: \mathcal P \to \nu \mathcal P$ is an equivalence; see (\ref{NFPro}), which clearly induces an equivalence $\nu\hspace{-1pt}:\hspace{-1pt} K^b(\mathcal P) \hspace{-1pt}\to\hspace{-1pt} K^b(\nu \mathcal P)$. It remains to prove Statement (2). By definition, we obtain binatural $R$-linear isomorphisms $\beta_{_{\hspace{-1pt}P,X}}: {\rm Hom}_{\mathfrak{A}}(X,\nu P)\rightarrow D {\rm Hom}_{\mathfrak{A}}(P,X)$, for all $P\in \mathcal P$ and $X\in \mf A$.

Fix a complex $X^\cdt$ over $\mf A$ and a bounded complex $P^\pdt$ over $\mathcal P$. We may define an $R$-linear map
$\beta_{_{\hspace{-1.3pt}P^\pdt, X^\ydt}}: \Hom_{\hspace{.5pt}C(\mathfrak{A})}(X^\ydt,\nu P^\pdt) \to
D{\rm Hom}_{\hspace{.5pt}C(\mathfrak{A})}(P^\pdt, X^\ydt)$ by setting \vspace{-2pt}
$$\beta_{_{\hspace{-1.3pt}P^\pdt, X^\ydt}}(\xi^\sdt)(\zeta^\sdt)= {\textstyle\sum}_{i\in \mathbb{Z}}\,(-1)^i\beta_{_{\hspace{-1.3pt}P^i, X^i}}(\xi^i)(\zeta^i),\vspace{-2pt}
$$ for $\xi^\sdt: X^\cdt\to \nu P^\pdt$ and $\zeta^\sdt: P^\pdt\to X^\cdt$ in $C(\mathfrak A).$ Using the binaturality of $\beta_{_{\hspace{-1pt}P,X}}$, we see that $\beta_{_{\hspace{-1.3pt}P^\pdt, X^\cdt}}$
is binatural and $\beta_{_{\hspace{-1.3pt}P^\pdt,X^\cdt}}(\xi^\sdt)(\zeta^\sdt)=0$ if $\xi^\sdt$ or $\zeta^\sdt$ is null-homotopic. This induces binatural $R$-linear maps $\bar \beta_{_{\hspace{-1.3pt}P^\pdt, X^\ydt}}: {\rm Hom}_{K(\mathfrak{A})}(X^\ydt,\nu P^\pdt) \to D {\rm Hom}_{K(\mathfrak{A})}(P^\pdt, X^\ydt)$ such that $\bar\beta_{_{\hspace{-1.3pt}P^\pdt, X^\ydt}}(\bar{\xi}^\sdt)(\bar{\zeta}^\sdt)=\beta_{_{\hspace{-1.3pt}P^\pdt,X^\cdt}}(\xi^\sdt)(\zeta^\sdt),$ which we claim are isomorphisms.

{\sc Sublemma.} {\it If $X^\ydt$ is a bounded complex, then $\bar\beta_{P^\pdt,X^\ydt}$ is an isomorphism.}

Indeed, we start with the case where $w(X^\ydt)= 1$, say $X^\ydt$ concentrates at degree $0$. Suppose first that $w(P^\pdt)=1.$ If $P^\pdt$ concentrates at degree $0$, then $\bar\beta_{_{\hspace{-1.5pt}P^\pdt,X^\pdt}}$ can be identified with $\beta_{_{\hspace{-1.5pt}P\hspace{.4pt}^0, X^0}}$, which is an isomorphism. Otherwise,
$\bar\beta_{_{\hspace{-1.5pt}P^\pdt,X^\ydt}}$ is a zero isomorphism. \vspace{-1pt} Suppose now that $w(P^\pdt) = s >1.$ By Lemma \ref{truncation-exact}, $K(\mathcal{P})$ has an exact triangle
$\xymatrixcolsep{18pt}\xymatrix{\hspace{-3pt} Q^\pdt\ar[r] &P^\pdt \ar[r] &L^\pdt \ar[r] &Q^\pdt[1]}\vspace{-3pt}$ with $w(Q^\pdt)< s$ and $w(L^\ydt)=1.$ \vspace{-1pt} This yields an exact triangle
$\xymatrixcolsep{18pt}\xymatrix{\nu Q^\ydt\ar[r] &\nu P^\pdt \ar[r] &\nu L^\ydt \ar[r] &\nu Q^\ydt[1]}\vspace{-1.5pt}$ in $K(\mf A)$. Since $I_{\hspace{-1pt}R}$ is injective, we obtain a commutative diagram with exact rows \vspace{-2pt}
$$\xymatrixcolsep{20pt}\xymatrixrowsep{16pt}\xymatrix{
(X^\pdt,\nu L^\ydt[-1])\ar[d]^{\bar\beta_{_{\hspace{-1.5pt}L^\ydt[-1], X^\ydt}}} \ar[r] & (X^\ydt, \nu Q^\ydt) \ar[d]^{\bar\beta_{_{\hspace{-1.5pt}Q^\ydt, X^\ydt}}} \ar[r] & (X^\ydt,\nu P^\pdt)\ar[d]^{\bar\beta_{_{\hspace{-1.5pt}P^\pdt, X^\ydt}}} \ar[r] & (X^\ydt,\nu L^\ydt) \ar[d]^{\bar\beta_{_{\hspace{-1.5pt}L^\ydt, X^\ydt}}}\ar[r] & (X^\ydt, \nu Q^\ydt[1])\ar[d]^{\bar\beta_{_{\hspace{-1.5pt}Q^\ydt[1], X^\ydt}}}\\
D(L^\ydt[-1],X^\ydt)\ar[r] & D(Q^\ydt, X^\ydt) \ar[r]&D(P^\pdt, X^\ydt) \ar[r] & D(L^\ydt, X^\ydt) \ar[r] & D (Q^\ydt[1], X^\ydt),}
\vspace{-2pt}$$
where $(M^\ydt, N^\ydt)=\Hom_{K(\mathfrak{A})}(M^\ydt, N^\ydt).$ By the induction hypothesis on $w(P^\pdt)$, we see that $\bar\beta_{_{\hspace{-1.5pt}P^\pdt, X^\ydt}}$ is an $R$-linear isomorphism.\vspace{-3pt}
Consider next the case where $w(X^\ydt)=t>1.$  By Lemma \ref{truncation-exact}, $K(\mathfrak{A})$ has an exact triangle $\xymatrixcolsep{20pt}\xymatrix{\hspace{-3pt}Z^\ydt\ar[r] &X^\ydt \ar[r] &Y^\ydt \ar[r] &Z^\ydt[1],}\vspace{-1.5pt}$ where $w(Z^\ydt)<t$ and $w(Y^\ydt)=1.$ This yields a commutative diagram with exact rows \vspace{-2pt}
$$\xymatrixcolsep{20pt}\xymatrixrowsep{16pt}\xymatrix{
(Z^\ydt[-1],\nu P^\pdt) \ar[d]^{\bar\beta_{_{\hspace{-1.5pt}P^\pdt,Z^\ydt[-1]}}}\ar[r] & (Y^\ydt,\nu P^\pdt) \ar[d]^{\bar\beta_{_{\hspace{-1.5pt}P^\pdt,Y^\ydt}}} \ar[r] & (X^\ydt,\nu P^\pdt) \ar[d]^{\bar\beta_{_{\hspace{-1.5pt}P^\pdt,X^\ydt}}} \ar[r] &
(Z^\ydt,\nu P^\pdt) \ar[d]^{\bar\beta_{_{\hspace{-1.5pt}P^\pdt,Z^\ydt}}} \ar[r] & (Y^\ydt[1],\nu P^\pdt)\ar[d]^{\bar\beta_{_{\hspace{-1.5pt}P^\pdt,Y^\ydt[1]}}} \\
D(P^\pdt, Z^\ydt[-1])\ar[r] & D(P^\pdt, Y^\ydt)\ar[r] & D(P^\pdt, X^\ydt)  \ar[r] & D(P^\pdt, Z^\ydt)\ar[r] & D(P^\pdt, Y^\ydt[1]).\vspace{-3pt} }$$
By the induction hypothesis, $\bar\beta_{_{\hspace{-1.5pt}P^\pdt,X^\pdt}}$ is an isomorphism.
This proves the sublemma.

In general, assume that $P^i=0$ for $i \not\in [m, n]$, where $m<n$. Considering the brutal truncations $M^\ydt=\kappa_{\ge m}(X^\ydt)$ and $N^\ydt=\kappa_{\le n}(M^\ydt)$ with canonical morphisms $\mu^\sdt: M^\ydt \to  X^\ydt$ and $\pi^\sdt: M^\cdt\to N^\ydt,$ we obtain a commutative diagram \vspace{-3pt}
$$\xymatrixcolsep{45pt}\xymatrixrowsep{16pt}\xymatrix{
\Hom_{K(\mathfrak{A})}(X^\ydt,\nu P^\pdt) \ar[d]_{\bar\beta_{_{\hspace{-1.5pt}P^\pdt, X^\ydt}}} \ar[r]^{\Hom(\bar\mu^\ydt\hspace{-2pt},  \,\nu P^\pdt)} & \Hom_{K(\mathfrak{A})}(M^\ydt,\nu P^\pdt) \ar[d]^{\bar\beta_{_{\hspace{-1.5pt}P^\pdt\hspace{-2pt}, \, M^\ydt}}}  &
\Hom_{K(\mathfrak{A})}(N^\ydt,\nu P^\pdt) \ar[d]^{\bar\beta_{_{\hspace{-1.5pt}P^\pdt, N^\ydt}}}\ar[l]_{\Hom(\bar\pi^\ydt\hspace{-2pt},\, \nu P^\pdt)}\\
D\Hom_{K(\mathfrak{A})}(P^\pdt, X^\ydt) \ar[r]^{D\Hom(P^\pdt\hspace{-2pt}, \,\bar\mu^\ydt)} & D\Hom_{K(\mathfrak{A})}(P^\pdt, M^\ydt) & D\Hom_{K(\mathfrak{A})}(P^\pdt, N^\ydt).\ar[l]_{D\Hom(P^\pdt\hspace{-2pt}, \,\bar\pi^\ydt)}
}\vspace{-3pt}$$
Since $P^i=0$ for $i\not\in [m, n]$, it is not difficulty to see that the horizontal maps are $R$-linear isomorphisms. Since $N^\ydt$ is bounded, by the sublemma, $\bar\beta_{_{\hspace{-1.2pt}P^\pdt, N^\ydt}}$ is an isomorphism, and so are $\bar\beta_{_{\hspace{-1.2pt}P^\pdt\hspace{-2pt}, \, M^\ydt}}$ and $\bar\beta_{_{\hspace{-1.2pt}P^\pdt, X^\ydt}}$. This establishes our claim. Then, by Lemma \ref{Hom-iso}, we obtain a binatural $R$-linear isomorphism \vspace{-3pt}
$$\tilde \beta_{_{\hspace{-1.3pt}P^\pdt, X^\ydt}}=D(\mathbb{L}_{P^\pdt\hspace{-.6pt}, \hspace{.4pt} X^\ydt }^{-1})\circ
\bar \beta_{_{\hspace{-1.3pt}P^\pdt, X^\ydt}}\circ \mathbb{L}_{X^\ydt\hspace{-.6pt}, \hspace{.4pt}\nu P^\pdt}^{-1}: {\rm Hom}_{D(\mathfrak{A})}(X^\ydt,\nu P^\pdt) \to D {\rm Hom}_{D(\mathfrak{A})}(P^\pdt, X^\ydt).\vspace{-2pt}$$
The proof of the proposition is completed.

\medskip

\noindent{\sc Remark.} The isomorphism stated in Proposition \ref{Hotopy-fun}(2) is known for the bounded derived category of a finite dimensional algebra; see \cite[Page 350]{H}.

\medskip

We are ready to obtain a sufficient condition for the existence of an almost split triangle in the derived categories of an abelian category with a Nakayama functor.


\begin{Theo}\label{ART-general}

Let $\mf A$ be an abelian $R$-category with $\mathcal P$ a strictly additive subcategory of projective objects of $\mf A$ and $\nu\hspace{-1pt}: \hspace{-1.5pt} \mathcal P \hspace{-1pt}\to \hspace{-1pt} \mf A$ a Nakayama functor. If $P^\pdt\in K^b(\mathcal P)$ and $\nu P^\pdt\in K^b(\nu \mathcal P)$ are strongly indecomposable, then $D^b(\mf A)$ \vspace{-1.6pt} has an almost split triangle
$\hspace{-3pt}\xymatrixcolsep{18pt}\xymatrix{\nu P^\pdt[-1] \ar[r] & M^\ydt \ar[r] & P^\pdt \ar[r] & \nu P^\pdt,}\hspace{-3pt}\vspace{-2.5pt}$ which is also almost split in $D(\mf A)$.

\end{Theo}

\noindent{\it Proof.} Assume that $P^\pdt\in K^b(\mathcal P)$ and $\nu P^\pdt\in K^b(\nu \mathcal P)$ are strongly indecomposable. By Lemma \ref{Hom-iso}, $P^\pdt$ and $\nu P^\pdt$ are strongly indecomposable in $D^b(\mf A)$. In view of Proposition \ref{Hotopy-fun}(2), we obtain an isomorphism \vspace{-3pt}
$$\Phi: \Ext_{D(\mathfrak{A})}^1(-, \nu P^\pdt[-1])=\Hom_{D(\mathfrak{A})}(-, \nu P^\pdt)\to D\Hom_{D(\mathfrak{A})}(P^\pdt, -),
\vspace{-2pt}$$ which restricts to an isomorphism \vspace{-2.5pt}
$$\Psi: \Ext_{D^b(\mathfrak{A})}^1(-, \nu P^\pdt[-1])=\Hom_{D^b(\mathfrak{A})}(-, \nu P^\pdt)\to D\Hom_{D^b(\mathfrak{A})}(P^\pdt, -).
\vspace{-2pt}$$

Choose a non-zero $R$-linear form $\theta: \End_{D^b(\mathfrak{A})}(P^\pdt)\to I{\hspace{-1pt}_R}$, which vanishes on the radical of $\End_{D^b(\mathfrak{A})}(P^\pdt)$. Then,
$\theta=\Psi_{P^\pdt}(\delta^\sdt)$ for some $\delta^\sdt\in \Ext_{D^b(\mathfrak{A})}^1(P^\pdt, \nu P^\pdt[-1])$. Since $\theta$ is in the right $\End_{D^b(\mathfrak{A})}(P^\pdt)$-socle of $D\End_{D^b(\mathfrak{A})}(P^\pdt)$, by Theorem \ref{AZE}, $\delta^\sdt$ is an almost-zero extension in $D^b(\mathfrak{A})$. On the other hand, since $D^b(\mathfrak{A})$ is a full triangulated subcategory of $D(\mathfrak{A})$, we see that $\delta^\sdt\in \Ext_{D(\mathfrak{A})}^1(P^\pdt, \nu P^\pdt[-1])$ such that $\Phi_{P^\pdt}(\delta^\sdt)=\theta$, which lies in the right $\End_{D(\mathfrak{A})}(P^\pdt)$-socle of $D\End_{D(\mathfrak{A})}(P^\pdt)$. Hence, $\delta^\sdt$ is an almost-zero extension in $D(\mathfrak{A})$. \vspace{-1pt} Therefore, $\delta^\sdt$ defines an almost split triangle
$\xymatrixcolsep{18pt}\xymatrix{\hspace{-3pt} \nu P^\pdt[-1] \ar[r] & M^\ydt \ar[r] &P^\pdt \ar[r] & \nu P^\pdt}\vspace{-1pt}\hspace{-2pt}$ in $D^b(\mf A)$, which is also an almost split triangle in $D(\mf A)$. The proof of the theorem is completed.

%

\medskip

\noindent{\sc Example.} Let $A$ be an $R$-algebra. By Lemma \ref{Alg-nf}, there exists a Nakayama functor $\nu_{\hspace{-2pt}_A}: {\rm proj}\hspace{.5pt}A \to \Mod\hspace{.5pt}A$, and hence, Theorem \ref{ART-general} applies in $D(\Mod A)$.

\medskip

As an application of Theorem \ref{ART-general}, we shall describe some almost split triangles in the derived categories of all modules over a locally finite dimensional algebra.

\begin{Cor}

Let $\La=kQ/J$ be a locally finite dimensional algebra over a field $k$, where $Q$ is locally finite and $J$ is weakly admissible.

\begin{enumerate}[$(1)$]

\item If $P^\pdt\in K^b({\rm proj}\hspace{.4pt}\La)$ is indecomposable, \vspace{-1.5pt} then $D^b({\rm Mod}\hspace{.4pt}\La)$ has an almost split triangle
$\xymatrixcolsep{18pt}\xymatrix{\hspace{-2pt}\nu_{\hspace{-1.5pt}_{\it\Lambda}}\hspace{-.5pt}P^\pdt[-1] \ar[r] & M^\pdt \ar[r] &P^\pdt \ar[r] & \nu_{\hspace{-1.5pt}_{\it\Lambda}} \hspace{-.5pt} P^\pdt\hspace{-1pt},}\vspace{-1pt}\hspace{-2pt}$ which is almost split in $D(\ModLa)$.

\item If $I^\cdt\in K^b({\rm inj}\hspace{.4pt}\La)$ is indecomposable, \vspace{-1.5pt} then $D^b({\rm Mod}\hspace{.4pt}\La)$ has an almost split triangle $\xymatrixcolsep{18pt}\xymatrix{\hspace{-2pt}I^\cdt\ar[r] & M^\pdt \ar[r] & \nu^{\mbox{-}}_{\hspace{-1.5pt}_{\it\Lambda}}I^\cdt[1] \ar[r] & I^\pdt[1], \hspace{-2pt}}\vspace{-1pt}$ which is almost split in $D(\ModLa)$. \vspace{-2pt}

\end{enumerate} \end{Cor}

\noindent{\it Proof.} Since $A$ is locally finite dimensional, ${\rm proj}\hspace{.4pt}\La$ is Hom-finite; see \cite[(3.2)]{BHL}, and so is $K^b({\rm proj}\hspace{.4pt}\La)$. Since the idempotents in $D^b(\Mod \La)$ split; see \cite[Corollary A]{LeC}, so do the idempotents in $K^b({\rm proj}\hspace{.4pt}\La)$. Thus, $K^b({\rm proj}\hspace{.4pt}\La)$ is Krull-Schmidt; see \cite[(1.1)]{LNP}. Moreover, the Nakayama functor $\nu_{\hspace{-1.5pt}_\mathit\Lambda}: {\rm proj}\hspace{.4pt}\La\to \Mod\La$ induces an equivalence $\nu_{\hspace{-1.5pt}_\mathit\Lambda}: K^b({\rm proj}\hspace{.4pt}\La)\to K^b({\rm inj}\hspace{.4pt}\La)$ with a quasi-inverse $\nu^{\hspace{.4pt}\mbox{-}}_{\hspace{-1.5pt}_{\it\Lambda}}: K^b({\rm inj}\hspace{.4pt}\La)\to
K^b({\rm proj}\hspace{.4pt}\La);$ see (\ref{Hotopy-fun}). If $P^\pdt\in K^b({\rm proj}\hspace{.4pt}\La)$ is indecomposable, then so is $\nu_{\hspace{-1.5pt}_{\it\Lambda}} P^\pdt\in K^b({\rm inj}\hspace{.4pt}\La)$. By Theorem \ref{ART-general}, there exists an almost split triangle as stated in Statement (1). On the other hand, if $I^\pdt\in K^b({\rm inj}\hspace{.4pt}\La)$ is indecomposable, then so is $\nu^{\mbox{-}}_{\hspace{-1.5pt}_{\it\Lambda}}I^\pdt[1]\in K^b({\rm proj}\hspace{.4pt}\La)$. By Theorem \ref{ART-general}, there exists an almost split triangle as stated in Statement (2). The proof of the corollary is completed.

\medskip

Similarly, we can describe some almost split triangles in the derived categories of all modules over a reflexive noetherian algebra.

\begin{Theo}\label{ART-NA}

Let $A$ be a reflexive noetherian $R$-algebra. Consider a strongly indecomposable complex $M^\ydt\in D^b(\Mod A)$.

\vspace{-2pt}

\begin{enumerate}[$(1)$]

\item If $M^\ydt$ is a complex over ${\rm mod}^+\hspace{-3pt}A$, \vspace{-1pt} then $D^b(\Mod A)$ has an almost split triangle $\xymatrixcolsep{16pt}\xymatrix{\hspace{-2pt}\ar[r] N^\ydt & L^\ydt \ar[r] & M^\ydt \ar[r] & N^\ydt[1]}\hspace{-2pt}\vspace{-3pt}$ if and only if $M^\ydt$ has a finite projective resolution $P^\pdt$ over ${\rm proj}\hspace{.5pt} A;$ and in this case, $N^\ydt\cong \nu_{\hspace{-2pt}_A}\hspace{-1pt}P^\pdt[-1]$, a complex over ${\rm mod}^{\hspace{.3pt}-\hspace{-3.2pt}}A.$

\vspace{1pt}

\item If $M^\ydt$ is a complex over ${\rm mod}^{\hspace{.3pt}-\hspace{-3.2pt}}A$, \vspace{-1pt} then $D^b(\Mod A)$ has an almost split triangle $\hspace{-4pt} \xymatrixcolsep{16pt}\xymatrix{M^\ydt \ar[r] & L^\ydt \ar[r] & N^\pdt \ar[r] & M^\ydt[1]\hspace{-4pt}}$ \vspace{-3pt} if and only if $M^\ydt$ has a finite injective co-resolution $I^\pdt$ over ${\rm inj}\hspace{.5pt} A;$ and in this case, $N^\ydt\cong \nu_{\hspace{-2pt}_A}^{\mbox{\hspace{.5pt}-\hspace{-2pt}}}I^\ydt[1]$, a complex over ${\rm mod}^+\hspace{-3pt}A$.

\end{enumerate}

\end{Theo}

\noindent{\it Proof.} Since $_AA$ is $R$-reflexive, by Lemma \ref{alg-ref-mod}, we see that ${\rm proj}\hspace{.5pt}A$ is Hom-reflexive over $R$. Thus, by Proposition \ref{Hotopy-fun}(1), the Nakayama functor $\nu_{\hspace{-2pt}_A}: {\rm proj}\hspace{.5pt}A\to \ModLa$ induces an equivalence $\nu_{\hspace{-2pt}_A}: K^b({\rm proj}\hspace{.5pt}A)\to K^b({\rm inj}\hspace{.5pt}A)$, which has a quasi-inverse
$\nu_{\hspace{-2pt}_A}^{\hspace{.8pt}\mbox{-\hspace{-3pt}}}: K^b({\rm inj}\hspace{.5pt}A)\to K^b({\rm proj}\hspace{.5pt}A)$. In particular, a complex $P^\pdt\in K^b({\rm proj}\hspace{.5pt}A)$ is strongly indecomposable if and only if $\nu_{\hspace{-2pt}_A}P^\pdt\in K^b({\rm inj}\hspace{.5pt}A)$ is strongly indecomposable. Moreover, by Theorem \ref{Noe-Ref}, ${\rm mod}^+\hspace{-3pt}A$ is an abelian category with enough projective modules in ${\rm proj}\hspace{.5pt}A$ and ${\rm mod}^-\hspace{-3pt}A$ is an abelian category with enough injective modules in ${\rm inj}\hspace{.5pt}A$, and consequently, every bounded complex over ${\rm mod}^+\hspace{-3pt}A$ has a projective resolution over ${\rm proj}\hspace{.5pt} A$ and every bounded complex over ${\rm mod}^-\hspace{-3pt}A$ has an injective co-resolution over ${\rm inj}\hspace{.5pt}A$; see \cite[(7.5)]{BGKHME}. Now, Statements (1) and (2) follow immediately from Theorems \ref{ART-nec} and \ref{ART-general}. The proof of the theorem is completed.

\medskip

{\sc Example.} If $R$ is a product of noetherian complete local commutative rings, then every noetherian $R$-algebra is reflexive; see \cite[Section 5]{AUS}. 

\medskip

In case $R$ is noetherian complete local, the finiteness of the global dimension of a noetherian $R$-algebra is related to the existence of almost split triangles in its derived category.

\begin{Cor}\label{ART-fgd}

Let $A$ be a noetherian $R$-algebra, where $R$ is a product of commutative noetherian complete local rings.
The following statements are equivalent.

\vspace{-2pt}

\begin{enumerate}[$(1)$]

\item The global dimension of $A$ is finite.

\vspace{.6pt}

\item Every indecomposable complex in $D^b({\rm mod}^{+\hspace{-3pt}}A)$ \vspace{.5pt} is the ending term of an almost split triangle in $D^b(\Mod A).$

\item Every indecomposable object in $D^b({\rm mod}^{-\hspace{-3pt}}A)$ is the starting term of an almost split triangle in $D^b(\Mod A).$

\end{enumerate} \end{Cor}

\vspace{-2pt}

\noindent{\it Proof.} First of all, $A$ is a reflexive noetherian $R$-algebra. Moreover, ${\rm mod}^{+\hspace{-3pt}}A$ is a Krull-Schmidt abelian subcategory of ${\rm RMod}\hspace{.5pt}A$; see \cite[Section 5]{AUS}, and by Theorem \ref{Noe-Ref}(2), so is ${\rm mod}^{-\hspace{-3pt}}A$. Thus, $D^b({\rm mod} \hspace{.5pt} A)$ and $D^b({\rm mod}^{\hspace{.5pt}\mbox{-}\hspace{-3.5pt}}A)$ are Krull-Schmidt; see \cite[Corollary B]{LeC}. Since ${\rm mod}^{+\hspace{-3pt}}A$ has enough projective modules in ${\rm proj}\hspace{.4pt}A$ and ${\rm mod}^{-\hspace{-3pt}}A$ has enough injective modules in ${\rm inj}\hspace{.4pt}A$, we see that $D^b({\rm mod} \hspace{.5pt} A)$ and $D^b({\rm mod}^{\hspace{.5pt}\mbox{-}\hspace{-3.5pt}}A)$ are full triangulated subcategories of $D(\Mod A)$; see \cite[(1.11)]{BaL}.

Let $A$ be of finite global dimension. Then, every bounded complex over ${\rm mod}^{+\hspace{-3pt}}A$ has a finite projective resolution over ${\rm proj}\hspace{.5pt}A$ and every bounded complex over ${\rm mod}^{\hspace{.5pt}\mbox{-}\hspace{-3.5pt}}A$ has a finite injective co-resolution over ${\rm inj}\hspace{.5pt}A$; see \cite[(7.5)]{BGKHME}. Thus, Statements (2) and (3) follow from Theorem \ref{ART-NA}. Conversely, assume that Statement (3) holds. Since ${\rm mod}^{-\hspace{-3pt}}A$ is Krull-Schmidt, we deduce from Theorem \ref{ART-NA}(2) that every module in ${\rm mod}^{-\hspace{-3pt}}A$ is of finite injective dimension. By Proposition \ref{Noe-Ref}(2), every module in ${\rm mod}^{+\hspace{-3pt}}A^{\rm op}$ is of finite projective dimension, and hence, $A^{\rm op}$ is of finite global dimension; see \cite[(9.12)]{ROT}. Being left and right noetherian as a ring, $A$ is of finite global dimension; see \cite[(9.23)]{ROT}. The proof of the corollary is completed.

\medskip

To conclude, we shall describe all possible almost split triangles in the bounded derived category of an abelian category with a Nakayama functor, enough projective objects and enough injective objects.

\begin{Theo}\label{ART-refl}

\hspace{-4pt} Let $\mf A$ be an abelian $R\hspace{-1pt}$-category with a Nakaya\-ma functor $\nu\hspace{-1pt}:\hspace{-1pt} \mathcal P \hspace{-2pt}\to \hspace{-1pt}\mf A$, where $\mathcal P$ is a Hom-reflexive strictly additive subcategory of projective objects of $\mf A$ such that $\mf A$ has enough projective objects in $\mathcal P$ and enough injective objects in $\hspace{-.5pt} \nu \mathcal P$.

\vspace{-1pt}

\begin{enumerate}[$(1)$]

\item If $Z^\ydt\in D^b(\mf A)$ is strongly indecomposable, \vspace{-1pt} then there exists an almost split triangle $\hspace{-2pt}\xymatrixcolsep{16pt}\xymatrix{X^\ydt \ar[r] & Y^\ydt \ar[r] & Z^\ydt \ar[r] & X^\ydt[1]}\hspace{-2pt}$  \vspace{-2pt} in $D^b(\mf A)$ if and only if $Z^\ydt$ has a finite projective resolution $P^\pdt$ over $\mathcal P;$ and in this case, $X^\ydt\cong \nu P^\pdt$.

\vspace{.5pt}

\item If $X^\ydt \hspace{-1.5pt} \in \hspace{-1.5pt} D^b(\mf A)$ is strongly indecomposable, \vspace{-1pt} then there exists an almost split triangle $\hspace{-3pt}\xymatrixcolsep{16pt}\xymatrix{X^\ydt \ar[r] & Y^\ydt \ar[r] & Z^\ydt \ar[r] & X^\ydt[1]}\hspace{-3pt}$ in $D^b(\mf A)$ \vspace{-2pt} if and only if $X^\ydt$ has a finite injective co-resolution $I^\ydt$ over $\nu \mathcal P;$ and in this case, $Z^\ydt\cong \nu^{\mbox{-}\hspace{-0.5pt}} I^\ydt$.

\vspace{.5pt}

\item If every object in $\mf A$ has a finite projective resolution over $\mathcal P ($respectively, injective co-resolution over $\nu\mathcal P)$, then $D^b(\mf A)$ has almost split triangles on the right $($res\-pectively, left$);$ and the converse holds in case $\mf A$ is Krull-Schmidt.


\end{enumerate} \end{Theo}

\noindent{\it Proof.} By Lemma \ref{Hotopy-fun}(1), the Nakayama functor $\nu: \mathcal P \to \mf A$ induces an equivalence
$\nu: K^b(\mathcal P) \to K^b(\nu \mathcal P)$ with a quasi-inverse $\nu^{\hspace{.4pt}\mbox{-}}: K^b(\nu \mathcal{P})\to K^b(\mathcal{P}).$ Since $\mf A$ has enough projective objects in $\mathcal P$ and  enough injective objects in $\nu \mathcal P$, every bounded complex over $\mf A$ has a projective resolution over $\mathcal P$ and an injective co-resolution over $\nu \mathcal P$; see \cite[(7.5)]{BGKHME}. In view of Theorems \ref{ART-nec} and \ref{ART-general}, we see easily that the first two statements hold true.

Next, assume that every object in $\mf A$ has a finite projective resolution over $\mathcal P$. Then, every bounded complex over $\mf A$ has a finite resolution over $\mathcal P$; see \cite[(7.5)]{BGKHME}. By Statement (1), $D^b(\mf A)$ has almost split triangles on the right. Conversely, suppose that $\mf A$ has almost split triangles on the right. In particular, every strongly indecomposable object in $\mf A$ is the ending term of an almost split triangle in $D^b(\mf A)$, and by Statement (1), it has a finite projective resolution over $\mathcal P$. If $\mf A$ in addition is Krull-Schmidt, then every object in $\mf A$ has a finite projective resolution over $\mathcal P$. This proves the first part of Statement (3), and the second part follows dually. The proof of the theorem is completed.

\medskip

\noindent{\sc Example.} (1) Let $A$ be an artin algebra over a commutative artinian ring $R$. Then ${\rm mod}\hspace{.5pt}A$ is a Hom-finite abelian $R$-category with enough projective modules and enough injective modules. Considering the Nakayama functor $\nu_{\hspace{-2pt}_A}: {\rm proj}\hspace{.5pt}A \to {\rm mod}\hspace{.5pt}A$, we see that Theorem \ref{ART-refl} applies in $D^b({\rm mod}\hspace{.5pt}A)$, and in particular, it includes Happel's results stated in \cite{Hap1}.

\vspace{1pt}

(2) Let $\La=kQ/J$ be a strongly locally finite dimensional algebra over a field $k$, where $Q$ is locally finite and $J$ is locally admissible. Then ${\rm mod}^{\hspace{.3pt}b\hspace{-2.5pt}}\La$ is a Hom-finite abelian $k$-category with enough projective modules in ${\rm proj}\,\La$ and enough injective modules in ${\rm inj}\,\La$. Considering the Nakayama functor $\nu_{\hspace{-2pt}_{\it\Lambda}}: {\rm proj}\hspace{.5pt}\La \to {\rm mod}^{\hspace{.3pt}b\hspace{-2.5pt}}\La$, we see that Theorem \ref{ART-refl} applies in $D^b({\rm mod}^{\hspace{.3pt}b\hspace{-2.5pt}}\La)$. In case $Q$ has no infinite path, every module in ${\rm mod}^{\hspace{.3pt}b\hspace{-2.5pt}}\La$ has a finite projective dimension and a finite injective dimension, and by Theorem \ref{ART-refl}(3), $D^b({\rm mod}^{\hspace{.3pt}b\hspace{-2.5pt}}\La)$ has almost split triangles.

\vspace{1pt}

(3) Let $\La=kQ/J$, where $Q$ is a locally finite quiver and $J$ is the ideal in $kQ$ generated by the paths of length two. Then,
every module in ${\rm mod}^{\hspace{.3pt}b\hspace{-2.5pt}}\La$ is of finite projective dimension over ${\rm proj}\,\La$ if and only if $Q$ has no right infinite path, and every module in ${\rm mod}^{\hspace{.3pt}b\hspace{-2.5pt}}\La$ is of finite injective dimension over ${\rm inj}\,\La$ if and only if $Q$ has no left infinite path. By Theorem \ref{ART-refl}, $D^b({\rm mod}^{\hspace{.3pt}b\hspace{-2.5pt}}\La)$ has almost split triangles $($on the left, on the right$)$ if and only if $Q$ has no $($left, right$)$ infinite path.

\end{document}